\documentclass[11pt]{amsart}

\usepackage{amssymb,amsmath}
\usepackage{graphicx}
\usepackage{bbm}
\usepackage{a4wide}

\theoremstyle{definition}
\newtheorem{example}{Example}
\newtheorem{remark}{Remark}

\newcommand{\ts}{\hspace{0.5pt}}
\newcommand{\nts}{\hspace{-0.5pt}}

\newcommand{\CC}{\mathbb{C}\ts}
\newcommand{\RR}{\mathbb{R}\ts}
\newcommand{\ZZ}{{\ts \mathbb{Z}}}
\newcommand{\SSS}{\mathbb{S}}
\newcommand{\TT}{\mathbb{T}}
\newcommand{\NN}{\mathbb{N}}
\newcommand{\XX}{\mathbb{X}}
\newcommand{\YY}{\mathbb{Y}}

\newcommand{\cA}{\mathcal{A}}
\newcommand{\cB}{\mathcal{B}}
\newcommand{\cC}{\mathcal{C}}
\newcommand{\cD}{\mathcal{D}}
\newcommand{\cE}{\mathcal{E}}

\newcommand{\cF}{\mathcal{F}}
\newcommand{\cO}{\mathcal{O}}

\newcommand{\cS}{\mathcal{S}}
\newcommand{\cU}{\mathcal{U}}

\newcommand{\vG}{\varGamma}
\newcommand{\vL}{\varLambda}
\newcommand{\one}{\mathbbm{1}}
\newcommand{\ii}{\mathrm{i}}
\newcommand{\ee}{\mathrm{e}}

\newcommand{\dd}{\, \mathrm{d}}

\newcommand{\MEF}{\mathsf{mef}}

\newcommand{\exend}{\hfill $\Diamond$}

\DeclareMathOperator{\dens}{dens}
\DeclareMathOperator{\vol}{vol}
\DeclareMathOperator{\card}{card}

\newcommand{\Hmm}[1]{\leavevmode{\marginpar{\tiny%
$\hbox to 0mm{\hspace*{-0.5mm}$\leftarrow$\hss}%
\vcenter{\vrule depth 0.1mm height 0.1mm width \the\marginparwidth}%
\hbox to 0mm{\hss$\rightarrow$\hspace*{-0.5mm}}$\\\relax\raggedright
#1}}}

\begin{document}

\title{Spectral notions of aperiodic order}

\author{Michael Baake}
\address{Fakult\"at f\"ur Mathematik, Universit\"at Bielefeld,
    Postfach 100131, 33501 Bielefeld, Germany}
\email{mbaake@math.uni-bielefeld.de}

\author{Daniel Lenz}
\address{Mathematisches Institut, Friedrich-Schiller-Universit{\"a}t
  Jena, 07743 Jena, Germany }
\email{daniel.lenz@uni-jena.de}

\begin{abstract}
  Various spectral notions have been employed to grasp the structure
  of point sets, in particular non-periodic ones. In this article, we
  present them in a unified setting and explain the relations between
  them. For the sake of readability, we use Delone sets in Euclidean
  space as our main object class, and give generalisations in the form
  of further examples and remarks.
\end{abstract}

\maketitle

\section{Introduction}

After the discovery of quasicrystals by Shechtman in 1982, which was
only published two years later \cite{Danny}, many people realised that
our common understanding of what `long-range order' might mean, is
incomplete (to put it mildly). In particular, little is known in the
direction of a classification, which --- despite the effort of many
--- still is the situation to date.  One powerful tool for the
analysis of order phenomena is provided by Fourier analysis, as is
clear from the pioneering work of Meyer \cite{M2}. Moreover, it not
surprising that methods from physical diffraction theory, most notably
the diffraction spectrum and measure of a spatial structure, have been
adopted and developed.

{}From another mathematical perspective, taking into account proper
notions of equivalence (which are needed for any meaningful
classification attempt), a similar situation is well-known from
dynamical systems theory. Here, the spectrum defined by Koopman
\cite{Koop} and later developed by von Neumann \cite{vN} and
Halmos--von Neumann \cite{HvN} led to a complete classification of
ergodic dynamical systems with pure point spectrum up to (metric)
isomorphism.

It is an obvious question how these spectral notions are related, and
part of this article aims at a systematic comparison, building on the
progress of the last 15 years or so. Since this means that large parts
of the paper will have review character, our exposition will be
informal in style. In particular, there will be no formal
theorems. Instead, we discursively present relevant statements,
concepts and underlying ideas and refer to the original literature for
more details and formal proofs as well as for generalisations. We hope
that the general ideas and results transpire more naturally this way,
and that the general flavour of the development is transmitted, too.
\smallskip

The paper is organised as follows. After the introduction of some
notions from point sets and spectral theory in
Section~\ref{Section:Prelim}, we begin with the diffraction spectrum
of an \emph{individual} Delone set in
Section~\ref{Section:individual}.  This part is motivated by the
description of the physical process of (kinematic) diffraction, where
one considers a single solid in a particle beam (photons or neutrons,
say) in order to gain insight into its internal structure, and by the
general mathematical aspects of Delone sets highlighted in
\cite{Lag-Delone}.  Next, in Section~\ref{Section:Diffraction}, we
extend the view by forming a dynamical system out of a given Delone
set and by extending the notion of the individual diffraction to that
of a diffraction measure of an (ergodic) dynamical \emph{system}.  The
two pictures (diffraction of individual sets and diffraction of
dynamical systems) are equivalent when the dynamical system is
uniquely ergodic, but we will not only look at this case.

Then, in a third step, we look at the dynamical spectrum of a Delone
dynamical system in Section~\ref{Section:The-dynamical-spectrum}, and
how it is related to its diffraction spectrum in
Section~\ref{sec:connections}. Beyond the equivalence in the pure
point case, which has been known for a while and is discussed in
Section~\ref{Section:Pure-point-diffraction}, we also look into the
more general case of mixed spectra, at least for systems of finite
local complexity (Section~\ref{Section:Factor}).  In this case, the
entire dynamical spectrum can still be described by diffraction.
However, one might have to consider the diffraction of a whole
\emph{family of systems} that are constructed from factors.

We then turn to the maximal equicontinous factor in
Section~\ref{Section:MEF}. This factor stores information on
continuous eigenfunctions. It can be used to understand a hierarchy of
Meyer sets via dynamical systems. Continuous eigenfunctions also play
a role in diffraction theory in the investigation of the so-called
Bombieri--Taylor approach.  Finally, in Section~\ref{Section:qpf}, we
have a look at our theory if the Delone set is replaced by suitable
quasiperiodic functions.  We compute autocorrelation and diffraction
in this case and discuss how the arising dynamical hull can be seen as
the maximal equicontinuous factor of the hull of a Delone set.
Moreover, we discuss an important difference between the diffraction
of quasiperiodic functions and that of Delone sets.

Our article gives an introduction to a field which has seen tremendous
developments over the last two decades, with steadily increasing
activity. In our presentation of the underlying concepts and ideas of
proofs, we do not strive for maximal generality but rather concentrate
(most of) the discussion to Delone sets and present examples in
various places. We have also included some pointers to work in
progess, as well as to some open questions. Part of the material, such
as the ideas concerning an expansion of sets into eigenfunctions in
Sections~\ref{Section:Pure-point-diffraction} and the discussion of
(diffraction of) quasiperiodic functions in Section \ref{Section:qpf},
do not seem to have appeared in print before (even though they are
certainly known in the community).

\section{Preliminaries}\label{Section:Prelim}

Let us begin by recalling some basic notions tailored to our later
needs.  We do not aim at maximal generality here but will rather
mainly be working in Euclidean space $\RR^{d}$. Some extensions will
be mentioned in the form of remarks.

We start with discussing point sets, see \cite[Sec.~2.1]{TAO} and
references therein for further details. A set consisting of one point
is called a \emph{singleton set}, while countable unions of singleton
sets are referred to as \emph{point sets}. A point set $\vL\subset
\RR^{d}$ is called \emph{locally finite} if $K\cap\vL$ is a finite set
(or empty), for any compact $K\subset\RR^{d}$.  Next, $\vL$ is
\emph{discrete} if, for any $x\in\vL$, there is a radius $r>0$ such
that $\vL\cap B_{r} (x) = \{ x \}$, where $B_{r} (x)$ denotes the open
ball of radius $r$ around $x$. If one radius $r>0$ works for all
$x\in\vL$, our point set is called \emph{uniformly discrete}. Next,
$\vL$ is called \emph{relatively dense} if a compact $K\subset
\RR^{d}$ exists such that $K+\vL=\RR^{d}$, where $A+B := \{ a+b \mid
a\in A , \ts b\in B\}$ denotes the Minkowski sum of two sets. Clearly,
if $\vL$ is relatively dense, there is a radius $R>0$ such that we see
the condition satisfied with $K = \overline{B_R (0)}$.

A \emph{Delone set} in $\RR^{d}$ is a point set that is both uniformly
discrete and relatively dense, so it can be characterised by two radii
$r$ and $R$ in the above sense.  They are therefore also called
$(r,R)$-sets in the literature. Delone sets are mathematical models of
atomic positions in solids, which motivates their detailed study in
our context.

A point set $\vL\subset\RR^{d}$ is said to have \emph{finite local
  complexity} (FLC) with respect to translations if, for any compact
neighbourhood $K$ of $0$, the collection of \emph{$K$-clusters of
  $\vL$},
\[
     \{ K \cap (\vL - x) \mid x \in \vL \}
\]
is a finite set. Again, it suffices to consider closed $R$-balls
around $0$ for all $R>0$, and $\vL$ is an FLC set if and only if $\vL
- \vL$ is locally finite; compare \cite[Prop.~2.1]{TAO}. Note that
clusters (or $R$-patches in the case we use a ball) are always defined
around a point of $\vL$, so that the empty set is \emph{not} a cluster
in our sense.

\begin{figure}
\begin{center}
  \includegraphics[width=0.8\textwidth]{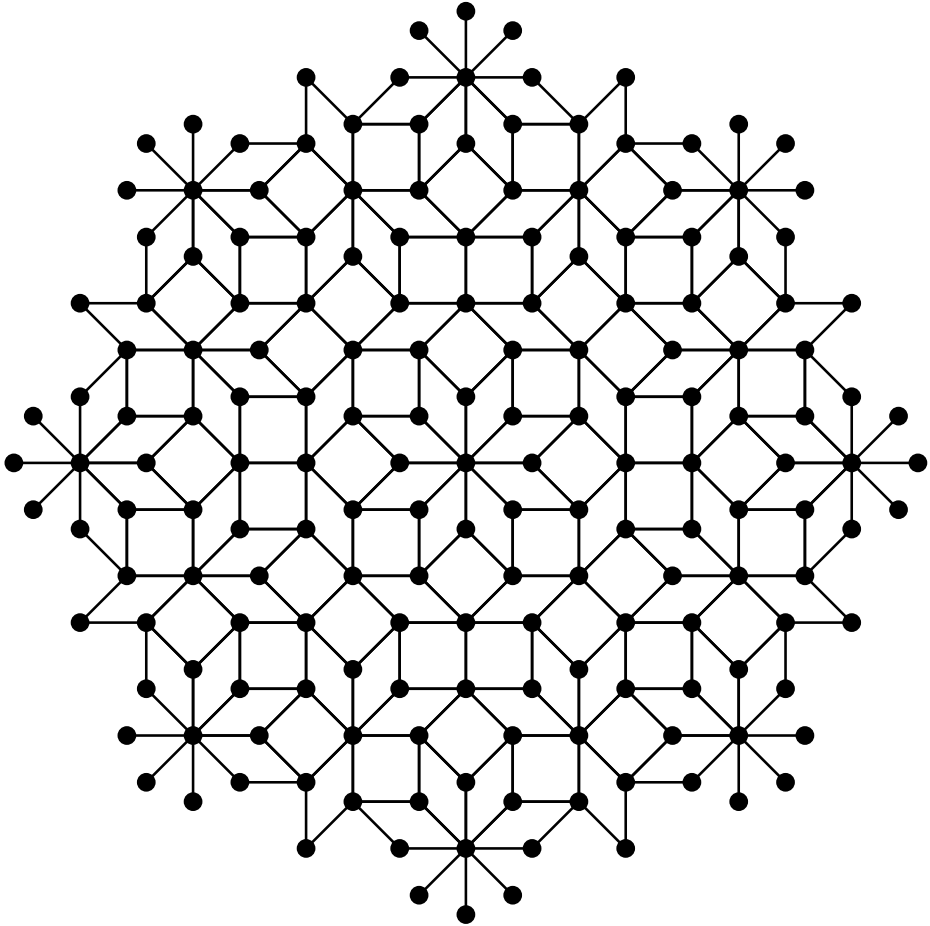}
\end{center}
\caption{A central patch of the eightfold symmetric Ammann--Beenker
  tiling, which can be generated by an inflation rule and it thus a
  self-similar tiling; see \cite[Sec.~6.1]{TAO} for details. The set
  of its vertex points is an example of a Meyer set, hence it is also
  an FLC Delone set. Moreover, it is a regular model set, as described
  in detail in \cite[Ex.~7.8]{TAO}.} \label{fig:ABtil} 
\end{figure}

A considerably stronger notion than that of a Delone set with finite
local complexity is that of a \emph{Meyer set}, where one demands that
$\vL$ is relatively dense while $\vL - \vL$ is uniformly discrete; see
\cite[Lemma~2.1 and Rem.~2.1]{TAO} for details and
\cite{M-Nato,M-beyond} for a thorough review. Clearly, every lattice
in Euclidean space is a Meyer set. Thus, Meyer sets can be thought of
as natural generalisations of lattices and this has been a very
fruitful point of view for the theory of Meyer sets. Meyer sets are
always subsets of model sets \cite{M2,M-Nato,TAO}, and important
idealisations of the atomic positions of quasicrystals.  A classic
example with eightfold symmetry in the plane is illustrated in
Figure~\ref{fig:ABtil}. \smallskip

There is a natural topology on the set of all Delone sets in Euclidean
space. This topology can be introduced in various ways. A very
structural way is to identify a Delone set $\vL$ with a measure by
considering its \emph{Dirac comb}
\[
     \delta^{}_{\! \vL} \, =  \sum_{x\in\vL} \delta^{}_{x} \ts ,
\]
where $\delta_{x}$ is the normalised point measure (or Dirac measure)
at $x$.  Clearly, different Delone sets correspond to different
measures. The vague topology on the measures then induces a topology
on the Delone sets \cite{BL-1}.  To identify Delone sets with measures
is more than a convenient mathematical trick. It is of great unifying
power as it allows us to treat sets, functions and measures on the
same footing. We will have more to say about this later.

At this point, we note that the topology on the Delone sets can be
generated by a metric as follows. Let
\[
    j \! : \, \SSS^d \xrightarrow{\quad}
          \RR^{d}\cup\{\infty\}
\]
be the stereographic projection. Here, $\SSS^d$ denotes the
$d$-dimensional sphere in $\RR^{d+1}$ and the point $\infty$ denotes
the additional point in the one-point compactification of $\RR^{d}$,
which is the image of the `north pole' under $j$. Let
$d^{}_{\mathrm{H}}$ be the Hausdorff metric on the set of compact
subsets of $\SSS^d$. Then, for any Delone set $\vL\subset \RR^{d}$,
the set $j^{-1} (\vL \cup\{\infty\}) $ is a closed and hence compact
subset of $\SSS^d$. Thus, via
\[
    d (\vL_1,\vL_2) \, := \,
      d^{}_{\mathrm{H}} \bigl( j^{-1} (\vL_1 \cup\{\infty\}),
          j^{-1} (\vL_2 \cup\{\infty\}) \bigr),
\]
we obtain a topology on the set of all Delone sets. It can be shown
that this is the same topology as the one discussed above \cite{LS}.
In this topology, the set of all $(r,R)$-Delone sets is compact
\cite{BL-1,LS}.

There is a canonical action of $\RR^{d}$ on the set of all Delone sets
by translations via
\[
    \RR^{d} \times \mbox{Delone sets}\xrightarrow{\quad}
    \mbox{Delone sets} \ts , \quad
     (t,\vL) \mapsto t + \vL \ts .
\]
Clearly, this action is continuous.

For any $(r,R)$-Delone set $\vL$,  its  \emph{hull}
\[
     \XX (\vL)  \, := \, \overline{\{ t + \vL \mid t \in \RR^{d} \}}
\]
is a closed and hence compact subset of the $(r,R)$-Delone sets.  By
construction, the hull is invariant under the translation action of
$\RR^{d}$. Thus, the pair consisting of the compact hull $\XX (\vL)$
and the restriction of the translation action of $\RR^{d}$ on this
hull is a \emph{dynamical system}, which we denote by $(\XX
(\vL),\RR^{d})$. As usual, this dynamical system is called
\emph{minimal} if the translation orbit of $\vL'$, which is $\{ t +
\vL' \mid t \in \RR^{d} \}$, is dense for every $\vL' \in \XX (\vL)$,
and it is called \emph{uniquely ergodic} if it possesses exactly one
probability measure which is invariant under the translation action.

The convolution $\varphi \ast \psi$ of $\varphi,\psi\in C_{\mathsf{c}}
(\RR^{d})$ is an element of $C_{\mathsf{c}} (\RR^{d})$ with
\[
    \bigl(\varphi \ast \psi\bigr) (x) \, :=
      \int_{\RR^{d}} \varphi (x-y) \, \psi (y) \dd y
\]
for all $x\in \RR^{d}$. We will identify measures on $\RR^{d}$ with
linear functionals on $C_{\mathsf{c}} (\RR^{d})$ by means of the
Riesz--Markov theorem. By the convolution of a measure $\nu$ with a
function $\varphi \in C_{\mathsf{c}} (\RR^{d})$, we mean the
continuous function $\nu * \varphi$ defined by
\[
    \bigl(\nu \ast \varphi) (x) \, =
     \int_{\RR^{d}} \varphi (x - y) \dd \nu (y) \ts .
\]
A particular role will be played by \emph{positive
definite measures}, which are measures $\nu$ with
\[
     \bigl(\nu \ast \widetilde{\varphi}\ast \varphi \bigr)(0)
     \, \geqslant \, 0
\]
for all $\varphi \in C_{\mathsf{c}} (\RR^{d})$, where
$\widetilde{\varphi}$ is defined by $\widetilde{\varphi} (x) =
\overline{\varphi (-x) }$. Any positive definite measure is
\emph{translation bounded}, meaning that $ \nu * \varphi$ is a bounded
function for all $\varphi \in C_{\mathsf{c}} (\RR^{d})$.

We will also need the Fourier transform of functions, measures and
distributions. For a complex-valued function $f$ on $\RR^{d}$ that is
integrable with respect to Lebesgue measure, we define its Fourier
transform $\widehat{f}$ as the complex-valued function given by
\[
   \widehat{f} (k) \, := \, \int_{\RR^{d}} \ee^{- 2 \pi \ii k x}
    f(x)  \dd x \ts .
\]
Clearly, this definition can be extended to finite measures; see
\cite[Ch.~8]{TAO} for details.  It turns out that it can also be
extended to various other classes of objects, including tempered
distributions. More delicate is the extension to unbounded measures,
where we refer to \cite{BF} for background. In particular, we note
that the Fourier transform of a positive definite measure exists and
is a positive measure.

\section{Diffraction of individual
objects}\label{Section:individual}

Here, we begin by considering a single Delone set $\vL\subset \RR^{d}$
and introduce and recall a spectral notion from the pioneering paper
\cite{Hof}, which is known as the \emph{diffraction measure} of $\vL$;
compare \cite[Sec.~9.1]{TAO} for a more detailed account.  In order to
put our approach in the general perspective of mathematical
diffraction theory, we will identify a Delone set $\vL$ with its Dirac
comb $\delta^{}_{\! \vL}$.  In our setting, the diffraction measure
emerges as the Fourier transform of the (natural) autocorrelation
measure, in extension of the classic Wiener diagram for integrable
functions; compare \cite[Sec.~9.1.2]{TAO}. Since $\delta^{}_{\! \vL}$
is an infinite measure, it cannot be convolved with itself, wherefore
one needs to proceed via restrictions to balls (or, more generally, to
elements of a general van Hove sequence \cite[Def.~2.9]{TAO}). Setting
$\delta^{R}_{\!\vL} := \delta^{}_{\!\vL \cap \overline{B^{}_{R}
    (0)}}$, we consider
\[
         \gamma^{R}_{\vL} \, := \, \frac{\delta^{R}_{\!\vL} *
         \widetilde{\delta^{R}_{\!\vL}}}{\vol (B^{}_{R} (0))}
\]
where $\widetilde{\mu}$ is the `flipped-over' version of a measure
$\mu$, defined by $\widetilde{\mu} (g) = \overline{\mu
  (\widetilde{g}\ts )}$ for $g\in C_{\mathsf{c}} (\RR^{d})$ and
$\widetilde{g}$ as above. Complex conjugation is not relevant in our
point set situation, but is needed for any extension to (complex)
weighted Dirac combs and general measures.

Every accumulation point of the family $\{ \gamma^{R}_{\vL} \mid R >
0\}$ in the vague topology, as $R\to\infty$, is called an
\emph{autocorrelation} of the Delone set $\vL$. By standard arguments,
compare \cite[Prop.~9.1]{TAO}, any Delone set possesses at least one
autocorrelation, and any autocorrelation is translation bounded. If
only one accumulation point exists, the autocorrelation measure
\[
       \gamma^{}_{\nts \vL} \, = \, \lim_{R\to\infty} \gamma^{R}_{\vL}
\]
is well-defined (we will only consider this situation later), and
called the \emph{natural autocorrelation}. Here, the term `natural'
refers to the use of balls as averaging objects, as they are closest
to the typical situation met in the physical process of
diffraction. In `nice' situations, the autocorrelation will not depend
on the choice of averaging sequences, as long as they are of van Hove
type (where, roughly stating, the surface to volume ratio vanishes in
the infinite volume limit). The volume averaged convolution in the
definition of $\gamma^{}_{\!\vL}$ is also called the \emph{Eberlein
  convolution} of $\delta^{}_{\!\vL}$ with its flipped over version,
written as
\[
     \gamma^{}_{\nts \vL} \, = \delta^{}_{\!\vL}\circledast
        \widetilde{\delta^{}_{\!\vL}} \ts .
\]
We refer to \cite[Sec.~8.8]{TAO} for some basic properties and
examples.

A particularly nice situation emerges when $\vL$ is an FLC set, so
$\vL-\vL$ is locally finite.  Then, assuming the natural
autocorrelation to exist, a short calculation shows that
\[
      \gamma^{}_{\nts \vL} \, = \sum_{z\in\vL-\vL} \! \eta(z) \, \delta_{z}
      \quad \text{with} \quad
      \eta (z) \, = \lim_{R\to\infty} \frac{\card \bigl( \vL^{}_{R}
          \cap (\vL^{}_{R} - z)\bigr) }{\vol{B^{}_{R} (0)}}.
\]
According to its definition, $\eta (z)$ can be seen as the frequency
of the vector $z$ from  the difference set $\vL - \vL$. Thus, the
autocorrelation of $\vL$ stores information on the set of difference
vectors of $\vL$ and their frequencies. Note that $\gamma^{}_{\nts \vL}$
in this case is a pure point measure on $\RR^{d}$.

By construction, the autocorrelation of any Delone set $\vL$ is a
positive measure, which is also positive definite.  As a consequence,
its Fourier transform, denoted by $\widehat{\gamma^{}_{\nts \vL}}$,
exists, and is a positive (and positive definite) measure. This
measure describes the outcome of a scattering (or diffraction)
experiment with our `idealised solid' when put into a coherent light
or particle source; see \cite{Cow} for background.  By continuity of
the Fourier transform, we have
\[
     \widehat{\gamma^{}_{\nts \vL}} \, =
     \lim_{R\to \infty} \widehat{ \gamma^{R}_{\vL}} \, =
     \lim_{R\to\infty} \, \frac{1}{\vol (B^{}_{R}(0))}
     \sum_{x,y\in \vL\cap B^{}_{R} (0)} e^{ 2 \pi \ii (x-y)(\cdot)}.
\]
Here, the function on the right hand side is considered as a measure
(namely the measure which has the function  as its density with respect
to Lebesgue measure)  and  the limit is taken in the sense of vague
convergence of measures.

Given the interpretation of the diffraction measure as outcome of a
diffraction experiment, it is natural that special attention is paid
to the set
\[
     \cB \, := \, \big\{ k\in\RR^{d} \mid
     \widehat{\gamma^{}_{\nts \vL}} \bigl( \{k\} \bigr) > 0 \big\} \ts .
\]
This set is denoted as the \emph{Bragg spectrum} (after the
fundamental contributions to structure analysis of crystals via
diffractions of the Braggs, father and son, which was honoured with
the Noble Prize in Physics in 1914).  The point measures of
$\widehat{\gamma^{}_{\nts \vL}}$ on the Bragg spectrum are known as
\emph{Bragg peaks}, and for any $k \in \cB$, the value
$\widehat{\gamma^{}_{\nts \vL}} (\{k\})$ is called the
\emph{intensity} of the Bragg peak.  In this context, there is the
idea around that one should have
\[
    \widehat{\gamma^{}_{\nts \vL} } (\{k\}) \, =
    \lim_{n\to \infty} \biggl| \frac{1}{ \vol (B^{}_{R} (0))}
    \sum_{x\in \vL\cap B^{}_{R} (0)} \ee^{2 \pi \ii k x} \biggr|^2 \ts .
\]
Indeed, this formula is quite reasonable as it says that the intensity
of the diffraction at $k$ is given as a square of a mean Fourier
coefficient. We will have more to say about its validity as we go
along.

\begin{remark}\label{Remark-BT}
  The validity of such a formula is discussed in \cite{Hof} with
  reference to work of Bombieri and Taylor \cite{BT1,BT2}, who used
  the formula without justification for certain systems coming from
  primitive substitutions. This was later justified in \cite{GK}. For
  regular model sets, the formula was shown in \cite{Martin}, but is
  also contained in \cite{M1}; see \cite[Prop.~9.9]{TAO} as well. In
  both cases, the special structure at hand is used. We will discuss a
  structural approach to it in Section \ref{Section:MEF}. More
  recently, the approach via amplitudes in the form of averaged
  exponential sums was extended to weak model sets of extremal
  density, where different methods have to be used; see
  \cite[Prop.~8]{BHS} for details. We shall come back to this topic
  later.  \exend
  \end{remark}

  From now on, whenever the meaning is unambiguous, we will drop the
  Delone set index and simply write $\gamma$ and $\widehat{\gamma}$
  for the autocorrelation and diffraction of $\vL$. Our approach is
  not restricted to Delone sets (see various remarks below), though we
  will mainly consider this case for ease of presentation.

\begin{example}\label{ex:integers}
  The set $\ZZ$ of integers, in our formulation, is described by the
  Dirac comb $\delta^{}_{\ZZ}$, and possesses the natural
  autocorrelation $\gamma = \delta^{}_{\ZZ}$, as follows from a
  straightforward Eberlein convolution; compare
  \cite[Ex.~8.10]{TAO}. Its Fourier transform is then given by
  $\widehat{\gamma} = \delta^{}_{\ZZ}$, as a consequence of the
  Poisson summation formula (PSF); see \cite[Ch.~9.2.2]{TAO} for
  details.

  More generally, given a crystallographic (or fully periodic) Delone
  set $\vL\subset\RR^d$, its Dirac comb is of the form
  $\delta^{}_{\!\vL} = \delta^{}_{S} * \delta^{}_{\vG}$ where $\vG=\{
  t\in\RR^{d} \mid t+\vG = \vG\}$ is the lattice of periods of $\vL$
  and $S$ is a finite point set that is obtained by the restriction of
  $\vL$ to a (true) fundamental domain of $\vG$; compare
  \cite[Prop.~3.1]{TAO}. Now, a simple calculation gives the natural
  autocorrelation
\[
      \gamma \, = \, \dens (\vG)
      \bigl( \delta^{}_{S} * \widetilde{\delta^{}_{S}} \ts
       \bigr) * \delta^{}_{\vG} \ts ,
\]
which is easily Fourier transformable by an application of the
convolution theorem together with the general PSF in the form
$\widehat{\delta^{}_{\vG}} = \dens (\vG) \, \delta^{}_{\vG^{*}}$,
where $\vG^{*}$ is the dual lattice of $\vG$. The result is the
diffraction measure
\[
     \widehat{\gamma} \, = \,
     \bigl( \dens (\vG) \bigr)^{2}  \lvert h
     \rvert ^{2} \delta^{}_{\vG^{*}}
\]
where $h = \widehat{\delta_{S}}$ is a bounded continuous function on
$\RR^{d}$; see \cite[Sec.~9.2.4]{TAO} for further details. We thus see
that the diffraction measure is a pure point measure that is
concentrated on the points of the dual lattice.

It is perhaps worth noting that the finite set $S$ in the above
decomposition of the Dirac comb $\delta^{}_{\!\vL}$ is not unique, and
neither is then the function $h$, because there are infinitely many
distinct possibilities to choose a fundamental domain of $\vG$.
Still, all functions $h$ that emerge this way share the property that
the values of $\lvert h \rvert^{2}$ agree on all points of $\vG^{*}$,
so that the formula for the diffraction measure is unique and
unambiguous.  \exend
\end{example}

\begin{remark}\label{rem:measure-dyn}
  As is quite obvious from our formulation, the Dirac comb of a Delone
  set is an example of a translation bounded measure on Euclidean
  space. This suggests that one can extend the entire setting to
  general translation bounded measures; compare \cite{BF,Hof,BL-1} as
  well as \cite[Chs.~8 and 9]{TAO}. Given such a measure, $\omega$
  say, one then defines its autocorrelation measure as
  $\gamma^{}_{\omega} = \omega \circledast \widetilde{\omega}$,
  provided this limit exists. It is then a translation bounded,
  positive definite measure, hence Fourier transformable by standard
  arguments \cite{BF}, and $\widehat{\gamma^{}_{\omega}}$ is a
  translation bounded, positive measure, called the \emph{diffraction
    measure} of $\omega$. This point of view was first developed in
  \cite{Hof}, and has been generalised in a number of articles; see
  \cite{TAO} and references therein for background, and \cite{LS} for
  a general formulation. \exend
\end{remark}

Figure~\ref{fig:ABspec} below shows an example of a diffraction
measure for an aperiodic point set, namely that of the Ammann--Beenker
point set introduced in Figure~\ref{fig:ABtil}. For the detailed
calculation in the context of regular cyclotomic model sets, we refer
to \cite[Secs.~7.3 and 9.4.2]{TAO}. \smallskip

Although the notion of a diffraction measure is motivated by the
physical process of diffraction, so that this approach looks very
natural for Delone sets as mathematical models of atomic positions in
a solid, the concept is by no means restricted to Delone sets, or even
to measures.

\begin{example}\label{ex:temp-distr}
  Let $\cS (\RR)$ denote the space of Schwartz functions on $\RR$ and
  $\cS' (\RR)$ its dual, the space of \emph{tempered distributions};
  see \cite{Schwartz,Wal} for general background. In this context,
  $\delta^{\prime}_{x}$ is a distribution with compact support,
  defined by $(\delta^{\prime}_{x} , \varphi) = - \varphi^{\ts\prime}
  (x)$, where we follow the widely used convention to write
  $(T,\varphi)$ for the evaluation of a distribution $T \in \cS'
  (\RR)$ at a test function $\varphi \in \cS (\RR)$. Note that
  $\delta^{\prime}_{x}$ is \emph{not} a measure. Tempered
  distributions of compact support are convolvable, and one checks
  that $\delta^{\ts\prime}_{x} * \delta^{\ts\prime}_{y} = \delta^{\ts
    \prime\prime}_{x+y}$.

  Let us now consider $\omega = \delta^{\prime}_{\ZZ} :=
  \sum_{x\in\ZZ} \delta^{\prime}_{x}$, which clearly is a tempered
  distribution. Also, we have $\omega = \delta^{\prime}_{0} *
  \delta^{}_{\ZZ}$, so that standard arguments imply the existence of
  the Eberlein convolution of $\omega$. A simple calculation gives
\[
     \gamma^{}_{\omega} \, = \, \omega \circledast \widetilde{\omega}
     \, = \, \delta^{\prime\prime}_{\ZZ} \, = \,
     \delta^{\prime\prime}_{0} *  \delta^{}_{\ZZ}\ts .
\]
This is a tempered distribution of positive type, so $(
\gamma^{}_{\omega}, \varphi * \widetilde{\varphi} \, ) \geqslant 0$
for all $\varphi \in \cS (\RR)$. Its Fourier transform, which always
exists as a tempered distribution, is then actually a positive
\emph{measure}, by an application of the Bochner--Schwartz theorem.
Observing that $\widehat{\ts\delta^{\prime\prime}_{0}\ts}$ is a
regular distribution, and thus represented by a smooth function, one
can check that
\[
     \widehat{\ts\delta^{\prime\prime}_{0}\ts} (y) \, = \, 4 \pi^{2} y^{2} .
\]
Now, using the convolution theorem together with the PSF, it is
routine to check that
\[
     \widehat{\gamma^{}_{\omega}} \, = \,
     \widehat{\delta^{\prime\prime}_{\ZZ}} \, = \,
     4 \pi^2 (.)^{2} \delta^{}_{\ZZ} \, =
     \sum_{y\in\ZZ} 4 \pi^2 y^2 \delta_{y} \ts .
\]
This is a positive pure point measure, the (natural) \emph{diffraction
  measure} of the tempered distribution $\omega$. In comparison to
previous examples, it is \emph{not} translation bounded, which makes
it an interesting extension of the measures in
Example~\ref{ex:integers}.

More generally, let us consider a lattice $\vG \subset \RR^{d}$. If $p
= (p^{}_{1}, \dots , p^{}_{d})$ denotes a multi-index (so all $p_{i}
\in \NN^{}_{0}$) with $\lvert p \rvert = p^{}_{1} + \ldots +\ts p^{}_{d}$
and $x^{p} = x^{p^{}_{1}}_{1} \cdots\ts x^{p^{}_{d}}_{d}$, as well as
the differential operator
\[
      D^{p} \, = \,\frac{\partial^{\lvert p \rvert}}
      {\partial x^{\ts p^{}_{1}}_{1} \cdots\ts \partial
        x^{\ts p^{}_{\nts d}}_{d}} \ts ,
\]
see \cite{Wal} for background, we get $\delta^{(p)}_{x} \! *
\delta^{(q)}_{y} = \delta^{(p+q)}_{x+y}$, where $(\delta^{(p)}_{x} ,
\varphi) := (-1)^{\lvert p \rvert} \bigl( D^{p} \varphi\bigr) (x)$ as
usual.  Now, for fixed $p$, consider the lattice-periodic tempered
distribution $\omega = \delta^{(p)}_{\vG} = \delta^{(p)}_{0} *
\delta^{}_{\vG}$. As before, the natural autocorrelation
$\gamma^{}_{\omega}$ exists, and is given by
\[
    \gamma^{}_{\omega} \, = \,
   \dens (\vG)\, \delta^{(2p)}_{\vG} \, = \, \dens (\vG) \,
   \delta^{(2p)}_{0} \nts * \delta^{}_{\vG} \ts .
\]
This is a tempered distribution of positive type again, so its Fourier
transform is a positive tempered measure.  Observing
\[
    \widehat{\delta^{(2p)}_{0}} (y) \, = \,
    (4 \pi^2)^{\lvert p \rvert} y^{2p}
\]
in analogy to above, one can employ the convolution theorem together
with the general PSF from Example~\ref{ex:integers} to calculate the
diffraction, which results in
\[
   \widehat{\gamma^{}_{\omega}} \, = \,
    ( 4 \pi^{2} ) ^{\lvert p \rvert} \,
   \dens (\vG)^{2} \, (.)^{2p} \, \delta^{}_{\vG^{*}}
   \, = \, \dens (\vG)^{2}
   \, (4 \pi^{2})^{\lvert p \rvert} \sum_{y\in \vG^{*}}
   y^{2p} \ts \delta^{}_{y} \ts .
\]
This measure is only translation bounded for $p=0$, where it reduces
to the diffraction measure of the lattice Dirac comb $\delta^{}_{\vG}$
of Example~\ref{ex:integers} as it must.

Due to the convolution structure, one can further generalise as
follows.  Let $\vL \subset \RR^{d}$ be a Delone set with natural
autocorrelation $\gamma^{}_{\!\vL}$, let $\nu$ be a tempered
distribution of compact support, and consider $\omega = \nu *
\delta^{}_{\!\vL}$. Clearly, this is a tempered distribution, with
existing (natural) autocorrelation.  The latter is given by
$\gamma^{}_{\omega} = (\nu * \widetilde{\nu}\,) * \gamma^{}_{\!\vL}$,
which is of positive type again. Fourier transform then results in the
diffraction
\[
     \widehat{\gamma^{}_{\omega}} \, = \,
     \lvert \widehat{\nu} \rvert^{2} \,
     \widehat{\ts\gamma^{}_{\!\vL}\ts}
\]
where $\widehat{\nu}$ is a smooth function on $\RR^{d}$.
\exend
\end{example}

\begin{remark}\label{rem:function-space}
  As one can see from the general structure of the volume-weighted
  convolution, the concept of a diffraction measure can be put to use
  in a wider context. Let us thus start from a locally convex space
  $\cF$ of functions on $\RR^{d}$ and let $\cF\ts '$ be its dual, the
  space of continuous linear functionals on $\cF$.  Examples include
  $C_{\mathsf c} (\RR^{d})$, which gives the regular Borel measures
  with the vague topology, and $\cS (\RR^{d})$, with the space $\cS\ts
  '(\RR^{d})$ of tempered distributions as its dual, but also the
  space $\cD (\RR^{d})$ of $C^{\infty}$-functions with compact
  support, then leading to the space $\cD\ts ' (\RR^{d})$ of
  distributions \cite{Schwartz,Wal}.  Various other combinations will
  work similarly.

  What we need is the concept of a functional of compact support, or a
  suitable variant of it, and the convolution of two linear
  functionals $G,H$ of that kind, as defined by
\[
      (G*H , \varphi) \, := \, (G \times H , \varphi^{\times})
\]
where $\varphi \in \cF$ and $\varphi^{\times} \! : \, \RR^{d} \times
\RR^{d} \xrightarrow{\quad} \CC$ is defined by $\varphi^{\times} (x,y)
= \varphi ( x+y)$. To expand on this, let us assume that a
distribution $F\in \cD\ts ' (\RR^{d})$ is given. Fix some $\varepsilon
> 0$ and let $c^{}_{r,\varepsilon} \in \cD (\RR^{d})$ be a
non-negative function that is $1$ on the ball $B_{r} (0)$ and $0$
outside the ball $\overline{B_{r+\varepsilon} (0)}$. Such functions
exist for any $r>0$. Now, consider
\[
     \gamma^{\, (r)}_{F,\varepsilon} \, := \, \frac{c^{}_{r,\varepsilon} F  *
     \widetilde{c^{}_{r,\varepsilon} F}}{\int_{\RR^{d}}
     c^{}_{r,\varepsilon} (x) \dd x}
\]
which is well-defined, with $\int_{\RR^{d}} c^{}_{r,\varepsilon} (x)
\dd x = \vol (B_{r} (0)) + \cO (1/r)$ as $r\to\infty$. If
$\lim_{r\to\infty} \gamma^{\, (r)}_{F,\varepsilon}$ exists and is also
independent of $\varepsilon$, which will be the case under some mild
assumptions on $F$, we call the limit $\gamma^{}_{F}$ the
\emph{natural autocorrelation} of the distribution $F$. More
generally, one can work with accumulation points as well. If $\gamma$
happens to be a tempered distribution, we are back in the situation
that $\widehat{\gamma^{}_{F}}$ is a positive measure, called the
\emph{diffraction measure} of the distribution $F$. This setting
provides a versatile generalisation of the diffraction theory of
translation bounded measures; see \cite{BLPS,ST} for a detailed
account. \exend
\end{remark}

\section{Diffraction of dynamical systems}\label{Section:Diffraction}

The diffraction measure of an individual Delone set is a concept that
emerges from the physical situation of a diffraction experiment. It is
both well founded and useful. Still, it has a number of shortcomings
that are related with the fact that it is not obvious how
$\widehat{\gamma}$ `behaves' when one changes the Delone set. Since
the mapping between autocorrelation and diffraction is Fourier
transform, and thus one-to-one, we can address this issue on the level
of the autocorrelation. Let us assume we have a Delone set $\vL$ whose
natural autocorrelation exists. Clearly, any translate of the set
should have the same autocorrelation, so
\[
      \gamma^{}_{t+\vL} \, = \, \gamma^{}_{\nts \vL}
      \quad \text{for all} \quad
      t\in\RR^{d} \ts ,
\]
and this is indeed a simple consequence of the van Hove property of
the family of balls $\{ B^{}_{R} (0) \mid R > 0\}$. In fact, a proof
only uses the (slightly weaker) F{\o}lner property of them for single
points.

Less obvious is what happens if one goes to the compact hull $\XX
(\vL)$ as introduced above.  Nevertheless, at least from a dynamical
systems point of view, it is very natural to define an autocorrelation
for a dynamical system. Here, one best starts with a measure-theoretic
dynamical system $(\XX (\vL), \RR^{d}, \mu)$, where $\mu$ is an
invariant probability measure on $\XX (\vL)$. In the large and
relevant subclass of uniquely ergodic Delone dynamical systems with
(FLC), the unique measure $\mu$ is the patch frequency measure.  Then,
any such measure-theoretic dynamical system $(\XX (\vL), \RR^{d},
\mu)$ comes with a autocorrelation $\gamma_{\mu}$ associated to it via
a closed formula (as opposed to a limit). This is discussed next,
where we follow \cite{BL-1}; see \cite{Gou-1} as well.

Choose a function $\chi \in C_{\mathsf{c}} (\RR^{d})$ and consider
the map $\gamma^{}_{\mu, \chi}\!  : \, C_{\mathsf{c}} (\RR^{d})
\xrightarrow{\quad} \CC $ defined by
 \[
       \varphi \,\mapsto \int_{\XX (\vL)}
       \sum_{x,y\in \vL'} \varphi (x - y) \, \chi (x)
       \dd \mu (\vL') \ts .
\]
Clearly, $\gamma^{}_{\mu,\chi}$ is a continuous functional on $
C_{\mathsf{c}} (\RR^{d})$. By the Riesz--Markov theorem, it can then
be viewed as a measure. Now, for fixed $\varphi \in C_{\mathsf{c}}
(\RR^{d})$, the map
\[
    C_{\mathsf{c}} (\RR^{d})\xrightarrow{\quad} \CC \ts ,
    \quad \chi \mapsto \gamma^{}_{\mu,\chi} (\varphi),
\]
is continuous. Hence, it is a measure as well. Moreover, as $\mu$ is
translation invariant, this measure can easily be seen to be invariant
under replacing $\chi$ by any of its translates. Hence, it must be a
multiple of Lebesgue measure.  Consequently, it will take the same
values for all $\chi$, which are \emph{normalised} in the sense that
they satisfy $\int_{\RR^{d}} \chi (t) \dd t = 1$. So, the map
$\gamma^{}_{\mu, \chi}$ will be independent of $\chi$ provided $\chi$
is normalised. Thus, we can unambiguously define
\[
     \gamma^{}_{\mu} \, := \, \gamma^{}_{\mu,\chi}
\]
for any such normalised $\chi$. This is then called the
\emph{autocorrelation} of the dynamical system $(\XX (\vL), \RR^{d},
\mu)$.

\begin{figure}
\begin{center}
  \includegraphics[width=0.86\textwidth]{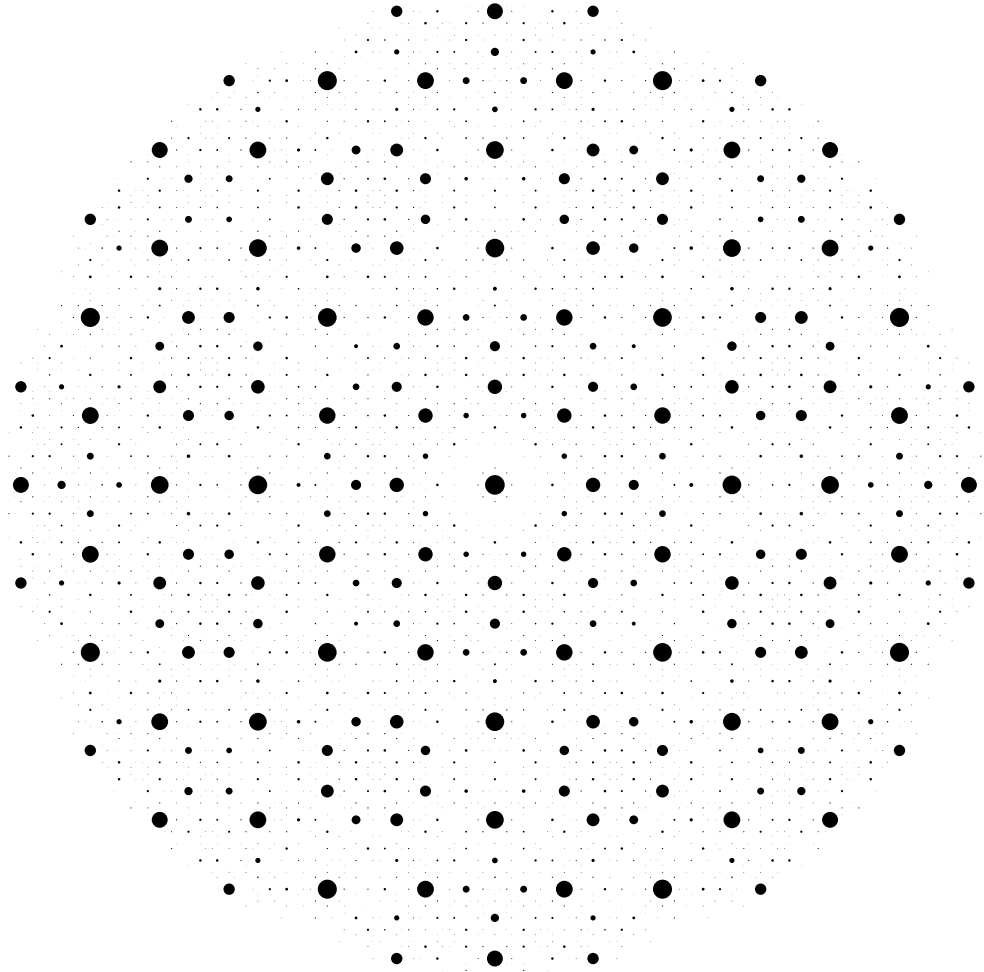}
\end{center}
\caption{Illustration of a central patch of the diffraction measure of
  the Ammann--Beenker point set of Figure~\ref{fig:ABtil}, which has
  pure point diffraction. A Bragg peak of intensity $I$ at $k\in\cB$
  is represented by a disc of an area proportional to $I$ and centred
  at $k$. Here, $\cB$ is a scaled version of $\ZZ[\ee^{\pi \ii/4}]$,
  which is a group; see \cite[Sec.~9.4.2]{TAO} for details. Clearly,
  $\cB$ is dense, while the figure only shows Bragg peaks beyond a
  certain threshold.  In particular, there are no extinctions in this
  case. At the same time, this measure is the diffraction measure of
  the Delone dynamical system defined by the (strictly ergodic) hull
  of the Ammann--Beenker point set, and $\cB$ is its dynamical
  spectrum.} \label{fig:ABspec}
\end{figure}

If $\mu$ is an ergodic measure, it can be shown that, for $\mu$-almost
every element $\vL'$ in the hull $\XX (\vL)$, the individual
autocorrelation $\gamma^{}_{\! \vL'}$ of $\vL'$ exists and equals
$\gamma^{}_{\mu}$. In general, the assessment of equality is
difficult, unless one knows that $\vL'$ is generic for $\mu$ in the
hull. However, if the dynamical system $(\XX(\vL),\RR^{d})$ is even
uniquely ergodic, the autocorrelation can be shown to exist and to be
equal to $\gamma^{}_{\mu}$ for every element in the hull. We refer to
\cite{BL-1} for further details and references. Important cases
include Delone sets derived from primitive substitution rules via
their geometric realisations, see \cite[Chs.~4 and 9]{TAO} for details
and many examples, and regular model sets in Euclidean space, such as
the Ammann--Beenker point set from Figures~\ref{fig:ABtil} and
\ref{fig:ABspec}; compare \cite[Chs.~7 and 9]{TAO} for more.

For any $(\XX (\vL),\RR^{d},\mu)$, the autocorrelation
$\gamma^{}_{\mu}$ can be shown to be a positive definite measure.
Hence, its Fourier transform exists and is a measure. This measure
will be called the \emph{diffraction measure} of the dynamical system,
and denoted by $\widehat{\ts\gamma^{}_{\nts \mu}\ts }$. As in the case
of the diffraction of an individual set, we will be particularly
interested in the point part of the diffraction measure.  The set of
atoms of this pure point part is again denoted by $\cB$ and called
\emph{Bragg spectrum}. It is then possible to compute the Bragg
spectrum via the following functions defined for each $k\in\RR^{d}$ by
\[
    c_{k}^{(R)} \! : \, \XX (\vL) \xrightarrow{\quad} \CC \ts ,
   \quad c_{k}^{(R)} (\vL'):= \frac{1}{ \vol (B^{}_{R} (0))}
    \sum_{x\in \vL'\cap B^{}_{R} (0)} \ee^{2 \pi \ii k x} .
\]
More specifically, as shown in \cite{Lenz},  we have
\[
    \widehat{\gamma^{}_{\mu}} ( \{k\}) \, = \lim_{R\to \infty}
    \|c_{k}^{(R)}\|^{2}_{L^2} \ts ,
\]
where $\|\cdot\|_{L^2}$ denotes the norm of the Hilbert space $L^2
(\XX (\vL), \mu)$,  and if  the dynamical system is ergodic, we even
have
\[
   \widehat{\gamma^{}_{\mu}} ( \{k\}) \,= \lim_{R\to \infty}
   \bigl| c_{k}^{(R)} (\vL') \bigr|^2
\]
for $\mu$-almost every $\vL'\in \XX (\vL)$. Note that, in these cases,
the corresponding limit will vanish for all $k\in\RR^{d} \setminus
\cB$.  One may expect that convergence holds for all $\vL' \in \XX
(\vL)$ in the uniquely ergodic case. However, this is not clear at
present.  We will have more to say about this in
Section~\ref{Section:MEF}.

\begin{remark}\label{rem:Choquet}
  In the preceding discussion, ergodicity of the measure on the hull
  has played some role. Thus, one may wonder about what happens for
  general measures. Thus, let $\nu$ be an arbitrary invariant
  probability measure on the hull that can be written as a convex
  combination $\nu = \sum_{i\in I} \alpha_{i} \, \mu_{i}$ of other
  invariant probability measures on the hull, hence $\alpha_{i} >0$
  and $\sum_{i\in I} \alpha_{i} = 1$. Then, using the same function
  $\chi$ for all autocorrelations, one sees that
\[
       \gamma^{}_{\nu} \, = \sum_{i\in I} \alpha_{i} \, \gamma^{}_{\mu_{i}} \ts .
\]
Invoking Choquet's theorem, compare \cite{Phelps} for background, one
can thus see that the analysis of the autocorrelations of extremal and
thus \emph{ergodic} invariant probability measures on the hull is the
essential step in the diffraction analysis of a Delone dynamical
system.  \exend
\end{remark}

\begin{remark}\label{rem:auto-relate} 
  At this point, we have discussed two ways of defining an
  autocorrelation, namely via a limiting procedure for individual
  Delone sets and via integration for hulls of Delone sets.  While
  these may seem very different procedures at first, we would like to
  stress that both have in common that they involve some form of
  \emph{averaging}. Indeed, in the limiting procedure, this is an
  average over $\RR^{d}$, while in the closed formula given above, it
  is an average over the hull. The connection between these two
  averages is then made by an ergodic theorem.  \exend
\end{remark}

\section{The dynamical spectrum}\label{Section:The-dynamical-spectrum}

In the preceding section, we have seen that any Delone dynamical
system $(\XX (\vL),\RR^{d},\mu)$ comes with an autocorrelation measure
$\gamma^{}_\mu$ (and thus also with a diffraction measure
$\widehat{\gamma^{}_\mu}$). We have also seen that this
autocorrelation measure agrees, for of a (typical) element of the
hull, with the individual autocorrelation of this element if the
measure $\mu$ is ergodic.  This suggests that there is a close
connection between properties of the dynamical system and the
diffraction. As was realised by Dworkin \cite{Dwo}, this is indeed the
case. This is discussed in this section. In order to discuss this
properly, we will first have to introduce the spectral theory of a
dynamical system. This is the spectral theory of what we call (in line
with various other people) the \emph{Koopman representation} of the
dynamical system, in recognition of Koopman's pioneering work
\cite{Koop}.  \smallskip

A Delone dynamical system $(\XX(\vL),\RR^{d},\mu)$ gives rise to a
unitary representation $T$ of $\RR^{d}$ on the Hilbert space $L^2
(\XX(\vL),\mu)$ via
\[
    T \! : \, \RR^{d}\xrightarrow{\quad}
    \mbox{unitary operators on $L^2 (\XX(\vL),\mu)$} \ts ,
    \quad t \mapsto T_t \ts ,
\]
with
\[
    T_t f \, = \, f(\cdot - t)\ts .
\]
Indeed, we obviously have $T_{t+s} = T_t T_s$ for any $t,s\in\RR^{d}$
as well as $T_0 = \one$. So, $T$ is a representation of $\RR^{d}$.
Also, as the measure $\mu$ is invariant, any $T_t$, with
$t\in\RR^{d}$, is isometric and, clearly, $T_{-t}$ is the inverse to
$T_t$. Thus, any $T_t$ is isometric and invertible and thus unitary.
Moreover, it is not hard to see that $T$ is strongly continuous, which
means that, for any fixed $f\in L^2 (\XX(\vL),\mu)$, the map
\[
   \RR^{d}\xrightarrow{\quad} L^2 (\XX(\vL),\mu)
   \ts , \quad t \mapsto T_t f \ts ,
\]
is continuous. We call the map $T$ the \emph{Koopman representation}
of the dynamical system.  As $T$ is a strongly continuous unitary
representation, Stone's theorem (compare \cite{Loomis}) guarantees the
existence of a projection-valued measure
\[
   E_T \! : \, \mbox{Borel sets on $\RR^{d}$}
   \xrightarrow{\quad} \mbox{projections
  on $L^2 (\XX(\vL),\mu)$}
\]
with
\[
     \langle f, T_t f \rangle \, =
      \int_{\RR^{d}} \ee^{2 \pi \ii t k}
      \dd  \rho^{}_f (k) \, = \, \widehat{\rho^{}_f} (-t) \ts ,
\]
for all $t\in\RR^{d}$,  where  $\rho^{}_f$ is the measure on $\RR^{d}$
defined by
\[
     \rho^{}_f (B) \, := \, \langle f, E_T (B)f \rangle\ts .
\]
The measure $\rho^{}_f$ is known as the \emph{spectral measure of $f$}
(with respect to $T$). It is the unique measure on $\RR^{d}$ with
$\langle f, T_t f \rangle = \widehat{\rho^{}_f} (-t)$ for all
$t\in\RR^{d}$. The study of the properties of the spectral measures is
then known as the \emph{spectral theory} of the dynamical system; see
\cite{Q} for a general exposition in the one-dimensional case.

In particular, the \emph{spectrum of the dynamical system} is
given as the support of $E$ defined by
\[
     \{ k\in\RR^{d} \mid E_T (B_\varepsilon (k)) \neq 0 \;\mbox{for all
      $\varepsilon >0$}\}.
\]
Of course, the spectrum is a set and as such does not carry any
information on the type of the spectral measures. For this reason, one
is mostly not interested in the spectrum alone, but also in
determining a spectral measure of maximal type (thus a spectral
measure having the same zero sets as $E$). We discuss a substitution
based system with mixed spectrum below in Example~\ref{ex:TM}. For us,
the following subset of the spectrum will be particularly
relevant. The \emph{point spectrum} of the dynamical system is given
as
\[
       \{ k\in\RR^{d} \mid E_T (\{k\}) \neq 0 \; \} \ts .
\]
A short consideration reveals that $k\in\RR^{d}$ belongs to the point
spectrum if and only if it is an eigenvalue of $T$. Here, an $f\neq 0$
with $f\in L^2 (\XX(\vL),\mu)$ is called an \emph{eigenfunction} to
the \emph{eigenvalue} $k\in\RR^{d}$ if
\[
    T_t f  \, = \, \ee^{2 \pi \ii t k} f
\]
holds for all $t\in \RR^{d}$.  Note that, following common practice,
we call $k$ (rather than $\ee^{2 \pi \ii k x}$) the eigenvalue, as
this matches nicely with the structure of the translation group as
well as its dual (the latter written additively).

If our dynamical system is ergodic, the modulus of any eigenfunction
must be constant (as it is an invariant function). So, in this case,
all eigenfunctions are bounded. If the system fails to be ergodic,
eigenfunctions need not be bounded. However, by suitable cut-off
procedures, one can always find bounded eigenfunction to each
eigenvalue; compare \cite{BL-1} for a recent discussion. It is not
hard to see that the eigenvalues form a group. Indeed,
\begin{itemize}
\item the constant function is an eigenfunction to eigenvalue $0$,
\item whenever $f$ is an eigenfunction to $k$, then $\overline{f}$ is
  an eigenfunction to $-k$, and
\item whenever $f$ and $g$ are bounded eigenfunctions to $k$ and
  $\ell$, respectively, the product $f g$ is an eigenfunction to $k + \ell$.
\end{itemize}
We denote this group of eigenvalues by $\cE (\mu)$. Standard reasoning
also shows that eigenfunctions to different eigenvalues are
orthogonal.  We will have more to say on eigenvalues and
eigenfunctions later.

\section{Connections between dynamical and diffraction
spectrum}\label{sec:connections}

Having introduced the dynamical spectrum, we now turn to the
connection with diffraction. The crucial ingredient is that the
Schwartz space $\cS(\RR^{d})$ can be embedded into $C (\XX (\vL))$ via
\[
    f \! : \, \cS(\RR^{d})\xrightarrow{\quad} C (\XX
   (\vL)) \ts , \quad \varphi \mapsto f_\varphi \ts ,
\]
with
\[
   f_\varphi (\vL') \, := \,\bigl( \varphi * \delta^{}_{\! \vL'}\bigr) (0)
   \, = \sum_{x\in \vL'} \varphi (-x) \, .
\]

\begin{remark}
  We could also work with the corresponding embedding of
  $C_{\mathsf{c}} (\RR^{d})$ into $C (\XX (\vL))$, and indeed this is
  often done. Note also that the existence of such embeddings will not
  be true for general dynamical systems, but rather requires the
  possibility of a `pairing' between the elements of the dynamical
  system and functions. Indeed, it is possible to extend (some of) the
  considerations below whenever such a pairing is possible
  \cite{BLPS,ST,LM}.  \exend
\end{remark}

Based on this embedding, one can provide the connection between
diffraction and dynamical spectrum. Here, we follow \cite{DM} (see
\cite{LM} as well), to which we refer for further details and proofs.
The key formula emphasised in \cite{Dwo} is
\[
   \bigl(\gamma^{}_\mu \ast \widetilde{\varphi}
   \ast \varphi\bigr) (0)   \, = \, \langle f_\varphi, f_\varphi \rangle
\]
for $\varphi \in \cS(\RR^{d})$.  This result was quite influential in
the field, as it highlighted a connection that was implicitly also
known in point process theory, compare \cite{Daley}, but had not been
observed in the diffraction context. Taking Fourier transforms and
using the denseness of $\cS(\RR^{d})$ in $L^2 (\RR^{d})$, one can use
this formula to obtain a (unique) isometric map
\[
    \Theta \! : \, L^2 (\RR^{d}, \widehat{\gamma^{}_{\mu}})
    \xrightarrow{\quad}  L^2 (\XX(\vL),\mu)\ts ,
    \quad \mbox{with} \; \Theta( \widehat{\varphi}\ts ) =
    f_\varphi
\]
for all $\varphi \in \cS(\RR^{d})$. Now, both $L^2$-spaces in question
admit a unitary representation of $\RR^{d}$. Indeed, we have already
met the Koopman respresentation $T$. Moreover, for any $t\in\RR^{d}$,
we have a unitary map
\[
    S_t \!  : \, L^2 (\RR^{d} , \widehat{\gamma^{}_{\mu}})
    \xrightarrow{\quad}
    L^2 (\RR^{d}, \widehat{\gamma^{}_{\mu}}) \ts ,
    \quad S_t h = \ee^{2 \pi \ii t (\cdot) } h,
\]
and these maps yield a representation $S$ of $\RR^{d}$ on the Hilbert
space $L^2 (\RR^{d},\widehat{\gamma^{}_\mu})$. Then, it is not hard to
see that $\Theta$ intertwines $S$ and $T$, which means that
\[
    \Theta S_{t} \, = \, T_{t}  \ts \Theta
\]
holds for all $t\in\RR^{d}$. In fact, this is clear when applying both
sides to functions of the form $\widehat{\varphi}$ for $\varphi \in
\cS (\RR^{d})$ and then follows by a denseness argument in the general
case. Consider now
\[
    \cU \, := \, \Theta \bigl( L^2
    (\RR^{d},\widehat{\gamma^{}_{\mu}}) \bigr) \nts \, = \,
   \overline{\mbox{Lin}\{ f_\varphi \mid \varphi \in \cS(\RR^{d})\}}
   \, \subset \, L^2 (\XX(\vL),\mu) \ts ,
\]
where the closure is taken in $L^2 (\XX(\vL),\mu)$. Then, $\cU$ is a
subspace. As $\Theta$ intertwines $S$ and $T$ and is an isometry, this
subspace is invariant under $T$ and the action of $S$ is equivalent to
the restriction of $T$ to this subspace.  In this sense, the
diffraction measure completely controls a subrepresentation of $T$.
\textbf{This is the fundamental connection between diffraction and
  dynamics.}

Using the map $\Theta$, we can easily provide a closed formula for the
pure point part of the diffraction measure.  Any $k\in \cB$ is an
eigenvalue of $S$ (with the characteristic function $1^{}_{\{k\}}$
being an eigenfunction). Hence, any $k\in \cB$ is an eigenvalue of $T$
with eigenfunction
\[
    c^{}_{k} \, := \, \Theta (1^{}_{\{k\}}) \ts .
\]
So, for any Bragg peak, there exists a canonical eigenfunction. This
is quite remarkable as eigenfunctions are usually only determined up
to some phase. The function $c^{}_{k}$ is not normalised in $L^2$.
Instead, using that $\Theta$ is an isometry, we obtain
\[
    \big\langle c^{}_{k}, c^{}_{k}\big\rangle^{}_{L^2 (\XX(\vL),\mu)} \, = \,
    \big\langle \Theta(1^{}_{\{k\}}), \Theta (1^{}_{\{k\}})
    \big\rangle^{}_{L^2 (\XX(\vL),\mu)} \, = \,
    \big\langle 1^{}_{\{k\}}, 1^{}_{\{k\}}
    \big\rangle^{}_{L^2 (\RR^{d},\widehat{\gamma^{}_{\mu}})}
    \, = \,  \widehat{\gamma^{}_{\mu}} ( \{k\} ) \ts.
\]
For the pure point part of the diffraction, we thus get
\label{diffraction-formula}
\[
    \bigl(\widehat{\gamma^{}_{\mu}}\bigr)_{\mathsf{pp}}
    \, = \sum_{k\in \cB} \|c^{}_{k}\|^2 \, \delta^{}_{k} \ts .
\]

For a given Delone set $\vL$, we have now considered two procedures to
investigate the associated diffraction, one via a limiting procedure
and one via considering the hull.  A short summary on how these two
compare may be given as follows:
\begin{eqnarray*}
\mbox{\textbf{point set $\vL$}} & \; \longleftrightarrow \;
& \mbox{\textbf{dynamical system $(\XX(\vL),\RR^{d},\mu)$}}\\
   \mbox{ $\gamma $ as a limit} &  \longleftrightarrow &
   \mbox{closed formula for $\gamma$}\\
   \mbox{$S$ on $L^2 (\RR^{d},\widehat{\gamma})$} &
    \longleftrightarrow & \mbox{restriction of $T$ to $\cU$}\\
   \mbox{Bragg spectrum}\; \cB & \xrightarrow{\quad\;}
&  \mbox{group of eigenvalues}\; \cE \\
    \mbox{Intensity}\; \widehat{\gamma}  (\{k\}) &
     \longleftrightarrow & \mbox{norm}\;  \|c^{}_{k} \|^2 .
\end{eqnarray*}
There is more to be said about the connection between the group
of eigenvalues and the Bragg spectrum, as we shall see later.

\section{Pure point diffraction and expansion in
eigenfunctions}\label{Section:Pure-point-diffraction}

The phenomenon of (pure) point diffraction lies at the heart of
aperiodic order, both in terms of physical experiments and in terms of
mathematical investigations. In this section, we take a closer look at
it.  \smallskip

We consider an ergodic Delone dynamical system $(\XX (\vL),
\RR^{d},\mu)$. This system comes with a unitary representation $T$ of
$\RR^{d}$ and a diffraction measure $\widehat{\gamma_\mu}$. It is said
to have \emph{pure point diffraction} if this measure is a pure point
measure.  It is said to have \emph{pure point dynamical spectrum} if
there exists an orthonormal basis of $L^2 (\XX(\vL),\mu)$ consisting
of eigenfunctions.

We have already seen in the previous section that the diffraction
measure controls a subspace of the whole $L^2 (\XX(\vL),\mu)$.
Accordingly, it should not come as a surprise that any $k\in \cB$ is
an eigenvalue of $T$ and pure point dynamical spectrum implies pure
point diffraction spectrum. Somewhat surprisingly it turns out that
the converse also hold. So, the Delone dynamical system $(\XX (\vL),
\RR^{d},\mu)$ has pure point diffraction if and only if it has pure
point dynamical spectrum.  Thus, the two notion of pure pointedness
are equivalent.

Following \cite{BL-1}, we can sketch a proof as follows: The
diffraction is pure point if $\widehat{\gamma}$ is a pure point
measure. By the discussion above, this is the case if and only if the
subrepresentation of $T$ coming from restricting to $\cU$ has pure
point spectrum. Clearly, if $T$ has pure point spectrum, then this
must be true of any subrepresentation as well and pure point
diffraction follows. To show the converse, note hat pure point
diffraction implies that all spectral measures $\varrho^{}_{\nts
  f_\varphi}$, with $\varphi \in \cS(\RR^{d})$, are pure point
measures (as these are equivalent to the spectral measures of
$\widehat{\varphi}$ with respect to $S$). We have to show that then
all spectral measures $\varrho^{}_{\nts f}$, with $f\in L^2
(\XX(\vL),\mu)$, are pure point measures. Consider
\[
    \cA \, := \, \{f\in C(\XX(\vL)) \mid
      \varrho^{}_f \;\mbox{is a pure point measure}\} \ts .
\]
Then, $\cA$ is a vector space with the following properties.
\begin{itemize}
\item It is an algebra. (This ultimately follows as the product of
  eigenfunctions is an eigenfunctions.)
\item It is closed under complex conjugation. (This ultimately follows
  as the complex conjugate of an eigenfunction is an eigenfunction.)
\item It contains all constant functions (as these are continuous
  eigenfunctions to the eigenvalue $0$).
\item It contains all functions of the form $f_\varphi$ (as has just
  been discussed) and these functions clearly separate the points of
  $\XX (\vL)$.
\end{itemize}
Given these properties of $\cA$, we can apply the Stone--Weierstrass
theorem to conclude that $\cA$ is dense in $C(\XX(\vL))$ with respect
to the supremum norm. Hence, $\cA$ is also dense in $L^2 (\XX
(\vL),\mu)$ with respect to the Hilbert space norm, and the desired
statement follows.

A closer inspection of the proof also shows that the group $\cE(\mu)$
of eigenvalues is generated by the Bragg spectrum $\cB$ if the system
has pure point diffraction spectrum.  Note that the Bragg spectrum
itself need not to be a group. The eigenvalues of $T$ which are not
Bragg peaks are called \emph{extinctions}. We refer to
\cite[Rem.~9.10]{TAO} for an explicit example. However, it is an
interesting observation in this context that $\cB$, in many examples,
actually \emph{is} a group, in which case one has identified the pure
point part of the dynamical spectrum as well.  This is the case for
the Ammann--Beenker point set, so that Figure~\ref{fig:ABspec} also
serves as an illustration of the dynamical spectrum.

In the case of pure point diffraction, the diffraction agrees with its
pure point part and the corresponding formula of the previous section
on p.~\pageref{diffraction-formula} gives $ \widehat{\gamma^{}_{\mu}}
= \sum_{k\in \cB} \|c^{}_{k}\|^2 \, \delta^{}_{k} $ with $c^{}_{k} =
\Theta (1^{}_{\{k\}})$. This way, the diffraction measure can actually
be used very efficiently to calculate the dynamical spectrum (in
additive formulation, as we use it here).

\begin{remark}
  The result on the equivalence of the two types of pure point
  spectrum has quite some history. As mentioned above, the work of
  Dworkin \cite{Dwo} provides the basic connection between the
  diffraction and the dynamical spectrum and gives in particular that
  pure point dynamical spectrum implies pure point diffraction
  spectrum; see \cite{Hof,Martin} for a discussion as well. In fact,
  for quite a while this was the main tool to show pure point
  diffraction spectrum \cite{Robbie,Sol}. For uniquely ergodic Delone
  dynamical systems with finite local complexity, the equivalence
  between the two notions of pure pointedness was then shown in
  \cite{LMS}. These considerations are modeled after a treatment of a
  related result for one-dimensional subshifts given in \cite{Q}.

  A different proof (sketched above), which permits a generalisation
  to arbitrary dynamical systems consisting of translation bounded
  measures, was then given in \cite{BL-1}. There, one can also find
  the statement that the Bragg spectrum generates the group of
  eigenvalues.  The setting of \cite{BL-1} does not require any form
  of ergodicity and applies to all Delone dynamical systems
  (irrespective of whether they are FLC or not), though it might be
  difficult then to actually determine the autocorrelation explicitly,
  despite the closed formula given in
  Section~\ref{Section:Diffraction}.

  A generalisation of \cite{LMS} to a large class of point processes
  was given in \cite{Gou-1}. This work applies to all Delone dynamical
  systems and requires neither ergodicity nor finite local
  complexity. In fact, it does not even require translation
  boundedness of the point process, but the weaker condition of
  existence of a second moment. A treatment containing both the
  setting of \cite{Gou-1} and \cite{BL-1} was then provided in
  \cite{LS} and, in a slightly different form, in \cite{LM}. These are
  the most general results up to date. The statement on the intensity
  of a Bragg peak being given by the square of an $L^2$-norm and the
  formula for $\widehat{\gamma^{}_{\mu}}$ can be found in \cite{Lenz}.

  It is worth noting that the equivalence between the dynamical
  spectrum and the diffraction spectrum only holds in the pure point
  case and does not extend to other spectral types, as follows from
  corresponding examples in \cite{vEM}; compare
  Section~\ref{Section:Factor} as well. It turns out, however, that
  --- under suitable assumptions --- the dynamical spectrum is
  equivalent to a \emph{family} of diffraction spectra
  \cite{BLvE}. Details will be discussed in
  Section~\ref{Section:Factor}. \exend
\end{remark}

We finish this section with a short discussion how pure point spectrum
can be thought of as providing an `Fourier expansion for the
underlying Delone sets'. To achieve this, we will need a normalised
version of the $c^{}_k$, with $k\in\cB$, given by
\[
    \widetilde{c^{}_k} \, := \, \frac{c^{}_{k}}
    {\bigl(\widehat{\gamma^{}_{\mu}} (\{k\})\bigr)^{1/2} }  \ts .
\]
As $\Theta$ is an isometry, we obtain from the very definition of
$c^{}_k$ for any $k\in\cB$ and any $\varphi \in \cS (\RR^{d})$
\[
    \big\langle f_\varphi \ts , \widetilde{c^{}_k}\big\rangle
    \, = \, \frac{ \big\langle \Theta(\widehat{\varphi}) \ts ,
     \Theta (1^{}_{\{k\}}) \big\rangle}
     {\bigl(\widehat{\gamma^{}_{\mu}} ( \{k\} )\bigr)^{1/2}}
     \, = \, \frac{ \big\langle \widehat{\varphi} \ts ,
               1^{}_{\{k\}} \big\rangle }
     {\bigl(\widehat{\gamma^{}_{\mu}} (\{k\}) \bigr)^{1/2}}
     \, = \, \bigl(\widehat{\gamma^{}_{\mu}} ( \{k\} )\bigr)^{1/2} \,
     \widehat{\varphi} (k).
\]
This in particular implies the relation
\begin{equation}\label{eq:fun-relation}
   \langle f_\varphi \ts , \widetilde{c^{}_k}\rangle
   \, \widetilde{c^{}_k} \, = \,
   \widehat{\varphi} (k) \, \widetilde{c^{}_k} \ts .
\end{equation}
This formula can be found in \cite{Lenz} (with a different proof).  It
will be used shortly.

Consider a Delone dynamical system $(\XX(\vL),\RR^{d},\mu)$ with pure
point diffraction and, hence, pure point dynamical spectrum.  The
basic aim is now to make sense out of the `naive' formula
\[
   \widehat{\delta^{}_{\!\vL'}} \, = \ts \sum_{k\in \cB}
   c^{}_{k} (\vL') \, \delta^{}_{k} \ts .
\]
To do so, we consider this equation in a weak sense. Thus, we pair
both sides with a $\widehat{\varphi}$ with some $\varphi \in \cS
(\RR^{d})$.  With $\varphi^{}_{\nts\text{-}}$ defined by
$\varphi^{}_{\nts\text{-}} (x) = \varphi (-x)$, we can then calculate
\[
\begin{split}
    \bigl(\widehat{\delta^{}_{\!\vL'}} \ts , \widehat{\varphi} \bigr)
    \; & = \; \bigl( \delta^{}_{\!\vL'} \ts , \varphi^{}_{\nts\text{-}} \bigr)
    \; = \; f^{}_{\varphi} (\vL') \; \overset{(!)}{=} \,
    \sum_{k\in \cB} \langle f^{}_{\varphi} \ts ,
    \widetilde{c^{}_{k}}\rangle \, \widetilde{c^{}_{k}} (\vL') \\ &
    \overset{\eqref{eq:fun-relation}}{=}
    \sum_{k\in\cB} \widehat{\varphi} (k) \, {c^{}_{k}} (\vL')
    \; = \,  \biggl( \, \sum_{x\in\cB} c^{}_{k} (\vL') \,
    \delta^{}_{k} \ts , \widehat{\varphi} \biggr) .
\end{split}
\]
Here, the second step follows from the definition of $f^{}_{\varphi}$,
while (!) requires some justification and this justification is
missing.  Note, however, that $(!)$ does hold in the $L^2$ sense.
Indeed, by the definition of pure point diffraction, the
$\widetilde{c^{}_{k}}$, with $k\in\cB$, form an orthonormal basis of
$\cU$ and hence
\[
    f^{}_{\varphi} \, = \ts \sum_{k \in \cB}
     \langle f_\varphi, \widetilde{c^{}_k}\rangle
     \,\widetilde{c^{}_k}
\]
is indeed valid. In fact, this is just the expansion of a function in
an orthonormal basis. So, the problem in the above reasoning is the
\emph{pointwise} evaluation of the Fourier series at $\vL'$.  We
consider this an intriguing open problem.

\begin{remark}
  Already in the original work of Meyer \cite{M1,M2}, it was an
  important point to capture the harmonic properties of point sets via
  trigonometric approximations. This led to the theory of harmonious
  sets; see \cite{M-Nato,M-beyond} for a detailed summary. More
  recently, Meyer has revisited the problem \cite{M3} and designed new
  schemes of almost periodicity that should help to come closer to a
  direct interpretation in the sense of an expansion.  \exend
\end{remark}

\section{Further relations between dynamical and diffraction
spectra}\label{Section:Factor}

Our approach so far was guided by the physical process of
diffraction. The latter is usually aimed at the determination (in our
terminology) of the Delone set, or as much as possible about it, from
the diffraction measure. This is a hard inverse problem, generally
without a unique solution. As mentioned before, diffraction is thus
tailored to one set, or to one dynamical system, and \emph{not}
invariant under (metric) isomorphism of dynamical systems.  This is
probably the reason why, from a mathematical perspective, it has not
received the attention it certainly deserves.

If one comes from dynamical systems theory, which has a huge body of
literature on spectral properties, it appears more natural to define a
spectrum in such a way that invariance under metric isomorphism is
automatic, and this was achieved by Koopman \cite{Koop}, and later
systematically explored by von Neumann \cite{vN}.  One celebrated
result in this context then is the Halmos--von Neumann theorem which
states that two ergodic dynamical systems with pure point spectrum are
(metrically) isomorphic if and only if they have the same spectrum,
and that any such system has a representative in the form of an
ergodic group addition on a locally compact Abelian group
\cite{vN,HvN, CFS}. This is well in line with the discussion of pure
point diffraction in the previous section (see \cite{LM} as
well). There, we have seen that pure point diffraction and pure point
dynamical spectrum are equivalent. So, in this case the diffraction
captures essentially the whole spectral theory.  A priori, it is not
clear what diffraction has to say on Delone dynamical systems with
mixed spectra, and the situation is indeed more complex then.

\begin{example}\label{ex:TM}
As was observed in \cite{vEM}, the subshift  $\XX^{}_{\mathrm{TM}}$
defined by the Thue--Morse (TM) substitution
\[
    \sigma^{}_{\mathrm{TM}} \! : \,
    a\mapsto ab \ts , \, b\mapsto ba \ts ,
\]
has a mixed dynamical spectrum that is \emph{not} captured by the
diffraction measure of the system; see also \cite[Secs.~4.6 and
10.1]{TAO} for a detailed discussion.  Note that we use a formulation
via substitutions here, but that one can easily obtain a Delone set as
well, for instance via using the positions of all letters of type $a$
in a bi-infinite TM sequence (letters correspond to unit intervals
this way).  To expand on the structure, the dynamical spectrum
consists of the pure point part $\ZZ \bigl[ \frac{1}{2}\bigr]$
together with a singular continous part that can be represented by a
spectral measure in Riesz product form,
\[
    \varrho^{}_{\mathrm{TM}} \, =
    \prod_{\ell=0}^{\infty} \bigl( 1 - \cos(2^{\ell+1} \pi x) \bigr) ,
\]
where convergence is understood in the vague topology (not pointwise)
and where $\varrho^{}_{\mathrm{TM}}$ is a spectral measure of maximal
type in the ortho-complement of the pure point sector.

Now, the diffraction measure picks up $\varrho^{}_{\mathrm{TM}}$
completely, but only the trivial part of the point spectrum, which is
$\ZZ$ in this case. Nevertheless, there is a single factor, the
so-called \emph{period doubling} subshift (as defined by the
substitution $\sigma^{}_{\mathrm{pd}} \! : \, a \mapsto ab \ts, \, b
\mapsto aa$), which has pure point spectrum (both diffraction and
dynamical).  Via the equivalence in this case, one picks up the entire
point spectrum, namely $\ZZ \bigl[ \frac{1}{2}\bigr]$. The period
doubling subshift emerges from the TM subshift via a simple sliding
block map; see \cite[Sec.~4.6]{TAO} for details.

In fact, it is possible to replace the TM system by a topologically
conjugate one, also based upon a primitive substitution rule (hence
locally equivalent in the sense of mutual local derivability), with
the property that one restores the equivalence of the two spectral
types for this system. The simplest such possibility emerges via the
induced substitution for legal words of length $2$; see
\cite[Sec.~5.4.1]{Q} or \cite[Sec.~4.8.3]{TAO} for details on this
construction. Here, this leads to a primitive substitution rule of
constant length over a $4$-letter alphabet. \exend
\end{example}

It turns out \cite{BLvE} that even in the case of mixed diffraction
one can capture the whole dynamical spectrum via diffraction (at least
in the case of systems with finite local complexity). However, one
will have to consider not only the diffraction of the original system
(which is not an isomorphism invariant) but also the diffraction of a
suitable set of factors (which when taken together provides an
isomorphism invariant). This is discussed in this section. We follow
\cite{BLvE}.

Let $(\XX(\vL),\RR^{d},\mu)$ be a Delone dynamical system of finite
local complexity. Let $T$ be the associated Koopman respresentation
and $E_T$ the corresponding projection valued measure.  A family $ \{
\sigma_\iota \}$ of measures on $\RR^{d}$ (with $\iota$ in some index
set $J$) is called a \emph{complete spectral invariant} when $E_T (A)
= 0$ holds for a Borel set $A\subset \RR^{d}$ if and only if
$\sigma_\iota (A) = 0$ holds for all $\iota \in J$.

An example for a complete spectral invariant is given by the family of
all spectral measures $\varrho^{}_f$, with $f\in L^2
(\XX(\vL),\mu)$. We will meet another spectral invariant shortly.
Recall that a dynamical system $(\YY,\RR^d)$ (i.e. a compact space
$\YY$ with a continuous action of $\RR^d$) is called a \emph{factor}
of $(\XX(\vL),\RR^d)$ if there exists a surjective continuous map
\[
     \Phi \! : \, \XX (\vL)\xrightarrow{\quad} \YY
\]
which intertwines the respective actions of $\RR^d$. In our context,
the dynamical systems will naturally be equipped with measures and we
will require additionally that the factor map maps the measure on $\XX
(\vL)$ onto the measure on $\YY$.

If $\YY $ is the hull of a Delone set with finite local complexity,
then $(\YY, \RR^{d},\nu)$ is called an FLC Delone factor. Of course,
any FLC Delone factor comes with an autocorrelation
$\gamma^{}_{(\YY,\RR^{d},\nu)}$ and a diffraction
$\widehat{\gamma}^{}_{(\YY,\RR^{d},\nu)}$. The main abstract result of
\cite{BLvE} then states that the family $
\widehat{\gamma}^{}_{(\YY,\RR^{d},\nu)}$, where $(\YY,\RR^{d},\nu)$
runs over all FLC Delone factors of $(\XX(\vL),\RR^{d},\mu)$, is a
complete spectral invariant for $T$.  In fact, it is not even
necessary to know the diffraction of all such factors. It suffices to
know the diffraction of so-called derived factors that arise as
follows.

Let $P$ be a $K$-cluster of $\vL$. For any $\vL' \in \XX (\vL)$, the
set of $K$-clusters of $\vL'$ is a subset of the $K$-clusters of
$\vL$, as a consequence of the construction of the hull $\XX (\vL)$.
We may thus define the \emph{locator set}
\[
   T^{}_{K,P} (\vL') \, = \, \{ t \in \RR^{d} \mid
    (\vL' - t)\cap K = P \}
   \, = \, \{ t \in \vL' \mid (\vL' - t)\cap K = P \}
   \, \subset \, \vL' \ts ,
\]
which contains the cluster reference points of all occurrences of $P$
in $\vL'$. Then, any $K$-cluster $P$ of $\vL$ gives rise to a factor
\[
   \YY  \, = \,  \YY_{K,P}
   \, := \, \{T^{}_{K,P} (\vL') \mid \vL' \in \XX (\vL) \}
\]
with factor map
\[
  \varPhi \, = \, \varPhi_{K, P}\nts : \; \XX
  \xrightarrow{\quad} \YY \ts , \quad
  X \mapsto T^{}_{K,P} (X) \ts .
\]
This factor will be called the factor \emph{derived from
  $(\XX,\RR^{d})$ via the $K\nts$-cluster $P$ of $\vL$}. It is the
diffraction of these factors (for all clusters) that is a complete
spectral invariant.

This result is relevant on many levels. On the abstract level, it
shows that the diffraction spectrum and dynamical spectrum are
equivalent in a certain sense. This may then be used to gather
information on the dynamical spectrum via diffraction methods. On the
concrete level, the result may even be relevant in suitably devised
experimental setups.

The considerations presented in this section raise naturally various
questions and problems. For example it seems that in concrete examples
often finitely many factors suffice. Thus, it would be of interest to
find criteria when this happens. Also, it is not unreasonable to
expect that in such situations also the diffraction of one factor (or
rather of one topologically conjugate system) suffices. Finally, it
would certainly of interest to extend the considerations to situations
where FLC does not hold.

\section{Continuous eigenfunctions and the maximal
equicontinuous factor} \label{Section:MEF}

Let $\vL$ be a Delone set with hull $\XX (\vL)$. Then, there is
natural embedding
\[
    \RR^{d}\xrightarrow{\quad} \XX (\vL) \ts ,
    \quad t\mapsto t + \vL,
\]
with dense range. In this way, the hull can be seen as a
compactification of $\RR^{d}$. As $\RR^{d}$ is an Abelian group, it is
then a natural question whether the hull carries a group structure
such that this natural embedding becomes a group homomorphism. In
general, this will not be the case.  Indeed, as shown in \cite{KL} for
an FLC Delone set $\vL$, such a group structure on the hull $\XX
(\vL)$ will exist if and only if $\vL$ is completely periodic. So, the
general question then becomes how close the hull is to being a
group. An equivalent formulation would be how much the metric on the
hull differs from being translation invariant. The concept of the
maximal equicontinuous factor (which we will recall below) allows one
to deal with these questions. This concept is not specific to Delone
dynamical systems. It can be defined for arbitrary dynamical systems
and this is how we will introduce it.

Throughout this section, we will assume that the occurring dynamical
systems are minimal (meaning that each orbit is dense). This is a
rather natural assumption as we want to compare the dynamical systems
to dynamical systems on groups, which are automatically minimal.
\smallskip

A dynamical system $(\TT, \RR^{d})$ is called a \emph{rotation on a
  compact group} if $\TT$ is a compact group and there exists an group
homomorphism
\[
    \xi \!  : \, \RR^{d} \xrightarrow{\quad} \TT
\]
with dense range  inducing the action of $\RR^{d}$ on $\TT$ via
\[
    t \cdot u \, := \, \xi (t) \ts  u
\]
for all $u\in\TT$ and $t\in\RR^{d}$. (Here, $\xi (t) u$ denotes the
product in the group $\TT$ of the two elements $\xi(t)$ and $u$.)  As
$\xi$ has dense range, the group $\TT$ must necessarily be Abelian. It
is well known (see e.g. \cite{ABKL} for a recent discussion) that any
rotation on a compact group is strictly ergodic (meaning uniquely
ergodic and minimal) and has pure point spectrum with only continuous
eigenfunctions (and the eigenvalues are just given by the dual of the
group $\TT$).

The \emph{maximal equicontinuous factor} (MEF) of a minimal dynamical
system $(\XX, \RR^{d})$ is then the largest rotation on a compact
group $(\TT,\RR^{d})$ which is a factor of $(\XX,\RR^{d})$.  It will
be denoted as $(\XX_{\MEF},\RR^{d})$ and the factor map will be
denoted as
\[
     \Psi_{\MEF} \! : \,
      \XX \xrightarrow{\quad} \XX_{\MEF} \ts .
\]
With this factor map at our disposal, the question of how close $\XX$
is to being a group becomes the question of how much $\Psi_{\MEF}$
differs from being bijective. In this context, one can naturally
distinguish three different regimes:

\begin{itemize}
\item The map $\Psi_{\MEF}$ is one-to-one everywhere (so, every point
  of $\XX_{\MEF}$ has exactly one inverse image). In this case, the
  hull carries the structure of a compact Abelian group (as it is
  isomorphic to $\XX_{\MEF}$).

\item The map $\Psi_{\MEF}$ is one-to-one almost everywhere
  (so, almost every point of $\XX_{\MEF}$ has exactly one inverse
  image). In this case, the hull is called an \emph{almost
    one-to-one extension} of its MEF.

\item The map $\Psi_{\MEF}$ is one-to-one in (at least) one point
  (meaning that there exists a point of $\XX_{\MEF}$ with exactly one
  inverse image). In this case, the hull is called an \emph{almost
    automorphic system}.
\end{itemize}

\begin{remark}
  Indeed, quite a substantial part of the general theory of the MEF is
  devoted to studying these three regimes \cite{Auslander88}. However,
  various other cases have been considered as well. This concerns in
  particular situations where the condition to be one-to-one is
  replaced by being $m$-to-one with a fixed integer $m$. In this
  context there is an emerging theory centered around the notion of
  coincidence rank; see \cite{ABKL} for a recent survey.  In the
  special case $m=2$, which occurs for instance for the TM subshift of
  Example~\ref{ex:TM} or for the twisted silver mean chain \cite{BG},
  interesting and strong results are possible because such an
  index-$2$ extension is quite restrictive; compare \cite{Hel} and
  \cite[Sec.~3.6]{Q} for background.  \exend
\end{remark}

Here, we are concerned with the situation that $\XX = \XX (\vL)$ is
the hull of a Delone set $\vL$. In this case, particular attention has
been paid to the case that $\vL$ is a Meyer set. In this case, the
corresponding parts of \cite{Aujogue,BLM,KS} can be summarised as
giving that these three regimes correspond exactly to the situation
that $\vL$ is crystallographic, a regular model set, a model set
respectively. We refrain from giving precise definitions or proofs but
rather refer the reader to \cite{ABKL} for a recent discussion; see
\cite{Kel} as well.

Next, we will provide an explicit description of the MEF for Delone
dynamical systems. In fact, it is not hard to see that a similar
description can be given for rather general dynamical systems as
well. For further details and reference we refer the reader to
\cite{ABKL}; see \cite{BLM} as well. Let $\cE_{\mathsf{top}}$ be the
set of continuous eigenvalues of $(\XX(\vL), \RR^{d})$. Here, an
eigenvalue $k\in \RR^{d}$ is called a \emph{continuous eigenvalue} of
$(\XX(\vL), \RR^{d})$ if there exists a continuous non-vanishing
function $f \! : \, \XX (\vL) \xrightarrow{\quad} \CC$ with
\[
      f (t + \vL') \, = \, \ee^{2 \pi \ii k t} f(\vL')
\]
for all $t\in\RR^{d}$ and $\vL'\in \XX (\vL)$. It is not hard to see
that the set of continuous eigenvalues is an (Abelian) group. We equip
this set with the discrete topology. Then, the Pontryagin dual
$\widehat{\cE_{\mathsf{top}}}$ of this group, which is the set
of all group homomorphisms
\[
    \cE_{\mathsf{top}} \xrightarrow{\quad}
    \SSS^{1} \, = \, \{z \in \CC : |z| =1\} \ts ,
\]
will be a compact group, In line with our previous convention, we
shall write this group additively and denote it by $\TT$.  There is a
natural group homomorphism
\[
    \xi \! : \, \RR^{d} \xrightarrow{\quad} \TT
    \quad \text{with } \, \xi(t) (k) := \ee^{2 \pi \ii t k}
\]
for all $t\in\RR^{d}$ and $k\in \cE_{\mathsf{top}}$. In this way,
$(\TT,\RR^{d})$ becomes a rotation on a compact Abelian group. Also,
$(\TT,\RR^{d})$ is a factor of $(\XX(\vL),\RR^{d})$. Indeed, choose
for each $k\in \cE_{\mathsf{top}}$ the unique continuous eigenfunction
$f_k$ with $f_k (\vL) =1$. Then, the map
\[
    \XX(\vL)\xrightarrow{\quad} \TT \, = \,
    \widehat{\cE_{\mathsf{top}}} \ts , \quad
   \vL'\mapsto (k\mapsto f_k (\vL')) \ts ,
\]
can easily be seen to be a factor map. Via this factor map, the
dynamical system $(\TT,\RR^{d})$ is the MEF of $(\XX(\vL),\RR^{d})$.

The preceding considerations show that there is a strong connection
between continuous eigenfunctions and the MEF. Somewhat loosely
speaking one may say that the MEF stores all information on continuous
eigenvalues.

In this context, dynamical systems coming from Meyer sets $\vL$ play a
special role. Indeed, this could already seen from the discussion
above relating a hierarchy of Meyer sets to injectivity properties of
the factor map $\Psi_{\MEF}$. It is also visible in recent results in
\cite{KS} showing that the dynamical system $(\XX (\vL),\RR^{d})$
coming from a Delone set with FLC has $d$ linearly independent
continuous eigenvalues if and only if it is conjugate to a dynamical
system $(\XX(\widetilde{\vL}),\RR^{d})$ with $\widetilde{\vL}$ a Meyer
set. In this sense, Delone dynamical systems with FLC and `many'
continuous eigenvalues are systems coming from Meyer sets.

Continuous eigenvalues also play a role in diffraction theory, as we
discuss next.  In Section~\ref{Section:Diffraction}, we have seen how
the autocorrelation of $(\XX(\vL),\RR^{d},\mu)$ can be computed by a
limiting procedure for $\mu$-almost every element $\vL' \in \XX (\vL)$
if $\mu$ is ergodic, and for all $\vL' \in \XX (\vL)$ if the system is
uniquely ergodic. In this context, we have also discussed the validity
of the formula
\[
    \widehat{\gamma  } (\{k\}) \, =
    \lim_{n\to \infty} \biggl| \frac{1}{ \vol (B^{}_{R} (0))}
    \sum_{x\in \vL'\cap B^{}_{R} (0)} \ee^{2 \pi \ii k x} \biggr|^2
\]
for almost every $\vL'$ in the ergodic case. Now, in the uniquely
ergodic case, this formula can be shown to hold even for all $\vL'$
provided  the eigenvalue $k$ is continuous \cite{Lenz}; see
\cite{Rob} for related earlier work as well.

\begin{remark}
  As discussed in Remark \ref{Remark-BT}, the validity of such a formula
  is known for sets coming from primitive substitutions as well as for
  regular model sets. In both cases, the associated Delone dynamical
  system is uniquely ergodic with only continuous eigenvalues. So, the
  mentioned work \cite{Lenz} provides a unified structural treatment.
  \exend
\end{remark}

It is an interesting open problem to which extent such a formula is
valid beyond the case of continuous eigenfunctions.  For example, it
is shown in \cite{Lenz} that such a formula holds for all linearly
repetitive systems even though such systems may have discontinuous
eigenfunctions \cite{BDM}. Also, the formula can be shown for weak
model sets of extremal density \cite[Prop.~6]{BHS}, where continuity
of eigenfunctions generally fails. It then also holds for generic
elements in the corresponding hull, equipped with a natural patch
frequency measure. Moreover, nonperiodic measures with locally finite
support \emph{and} spectrum, as recently constructed in \cite{M-new},
are further examples with well-defined amplitudes.  So, there is room
for generalisation, and hence work to be done to clarify the
situation.

\section{Quasicrystals and hulls of quasiperiodic
functions}\label{Section:qpf}

So far, we have (mostly) considered the dynamical system $(\XX (\vL),
\RR^{d})$ arising from a Delone set $\vL$. Special emphasis has been
paid to the case that this system is minimal and uniquely ergodic with
pure point point dynamical spectrum and only continuous
eigenfunctions. Indeed, these are the systems to which all results of
the preceding four sections apply.  In particular, these systems have
pure point diffraction and the set of (continuous) eigenvalues is a
group generated by the Bragg spectrum and
\[
     \Psi_{\MEF} \! : \, \XX (\vL)\xrightarrow{\quad} \TT
\]
is the factor map to its MEF, where $\TT$ is given as the dual group
of the group of eigenvalues. While it is not clear at present what a
mathematical definition for a quasicrystal should be, it seems
reasonable that such systems should fall into the class of
quasicrystals. At the same time, certain quasiperiodic functions are
also sometimes treated under the label of quasicrystals. In this
section, we compare these two approaches and also compute the
diffraction of a quasiperiodic function. This will actually show an
important structural difference in the diffraction measure which seems
to favour Delone sets as mathematical models for quasicrystals over a
description via quasiperiodic functions.  \bigskip

Let $\cC$ be a countable subset of $\RR^{d}$. Let $a^{}_{k}$, $k\in
\cC$, be non-vanishing complex numbers that satisfy the summability
condition
\[
   \sum_{k\in \cC } | a^{}_{k} | \, < \, \infty \ts .
\]
Denote the subgroup of $\RR^{d}$ generated by $\cC$ by $\cE'$, so
$\cE' = \langle \cC \rangle$. Define
\[
     u \! : \, \RR^{d}\xrightarrow{\quad} \CC \ts ,
     \quad  u (x) \, = \sum_{k\in \cC}
     a^{}_{k} \, \ee^{  2 \pi \ii k x} .
\]
By the summability condition, the sum is absolutely convergent and the
function is continuous and bounded. In fact, such functions are known
as \emph{quasiperiodic functions} in the sense of Bohr; see
\cite{Cord,Katz}, or \cite[Sec.~8.2]{TAO} for a short summary.

Clearly, we can view a bounded continuous function $f$ as a
Radon--Nikodym density relative to Lebesgue measure, and then identify
$f$ with the translation bounded measure defined that way. 
Consequentliy, we can equip the set of such functions with the
vague topology induced from measures. In particular, we can consider
the \emph{hull} of $f$ defined by
\[
    \XX (f) \, := \, \overline{\{ f(\cdot - t) \mid t\in\RR^{d}\} } \ts ,
\]
where the closure is taken in the vague topology on measures. Then,
$\XX (f)$ is compact and $\RR^{d}$ acts continuously via translations
on it (see e.g. \cite{BL-1}). Thus, we are given a dynamical system
$(\XX (f),\RR^{d})$. Assume now that $f = u$ is the quasiperiodic
function introduced above. Then, the closure $\XX (u)$ actually agrees
with the closure of the translates of $u$ in the topology of uniform
convergence, so
\[
    \XX (u) \, = \,  \overline{\{ u(\cdot - t) \mid
    t\in \RR^{d} \} }^{\|\cdot\|_\infty} .
\]
In particular, all elements in $\XX (u)$ (which are apriori only
measures) are continuous bounded functions.  Moreover, by standard
theory of almost periodic functions, compare \cite{Loomis} or
\cite[Sec.~8.2]{TAO} and references given there, this closure has the
structure of an Abelian group. More specifically, define
\[
    \xi \! : \, \RR^{d} \xrightarrow{\quad} \XX (u) \ts ,
    \quad t\mapsto u ( \cdot -t) \ts .
\]
Then, there exists a unique group structure on $\XX (u)$ making $\xi$
a homomorphism of Abelian groups (see \cite{LR} as well for a recent
discussion). This homomorphism has dense range and the translation
action of $\RR^{d}$ on $\XX (u)$ is given by
\[
   \RR^{d}\times \XX (u)\xrightarrow{\quad} \XX (u) \ts ,
   \quad (t, v) \mapsto v(\cdot - t)  = \xi (t) \cdot v \ts .
\]
Thus, $(\XX(u),\RR^{d})$ is a rotation on a compact group (in the
notation of Section~\ref{Section:MEF}). In particular, it is strictly
ergodic and has pure point dynamical spectrum with only continuous
eigenfunctions. Now, it is not hard to see that
\[
   \cC \, = \, \Big\{ k\in \RR^{d} \, \Big| \lim_{R\to \infty}
   \frac{ \int_{B^{}_{R} (0)}  u (x) \, \ee^{- 2 \pi \ii k x} \dd x }
   {\vol (B^{}_{R} (0))} \neq 0 \Big\}  .
\]
So, by standard theory of almost periodic functions, we infer
\[
   \cE' \, = \,  \langle \cC \rangle
    \, = \, \widehat{\XX (u)} \ts .
\]
Dualising once more we infer
\[
   \XX (u) \, = \, \widehat{\cE'} \ts .
\]

Assume now that the group $\cE'$ is the group of eigenvalues of the
uniquely ergodic minimal system $(\XX (\vL),\RR^{d})$ with pure point
spectrum, which has only continuous eigenfunctions. Then, its dual
group $\widehat{\cE'}$ is the MEF of $(\XX (\vL),\RR^{d})$, as
discussed at the beginning of this section. Moreover, as just derived,
this dual group is isomorphic to $\XX (u)$. Putting this together, we
see that the map $\Psi_{\MEF}$ can be considered as a map
\[
     \Psi_{\MEF} \! : \, \XX (\vL)\xrightarrow{\quad} \XX (u) \ts .
\]
In terms of the associated dynamical systems, we thus find a precise
relationship between the hulls of $\vL$ and of $u$: One is a factor of
the other and, in fact, a special one via the connection with the MEF.

These considerations can be slightly generalised as follows. Let $(\XX
(\vL),\RR^{d})$ be uniquely ergodic with pure point spectrum, and only
continuous eigenfunctions and group of eigenvalues $\cE$. If the group
$\cE' = \langle \cC \rangle$ is only a subgroup of $\cE$, we would
still get a factor map
\[
     \Psi \! : \, \XX(\vL)\xrightarrow{\quad} \XX (u) \ts ,
\]
as, in this case, the dual of the group $\cE'$ can easily be seen to
be a factor of $\widehat{\cE}$.

\begin{remark}
  The preceding considerations naturally raise the question whether,
  to any countable set $\cC$ and the induced group $\cE'$, one can
  find a uniquely ergodic minimal Delone dynamical system with pure
  point spectrum, only continuous eigenvalues and dynamical spectrum
  $\cE'$. The answer to this question is positive. In fact, it is even
  possible to find a Meyer set $\vL$ such that its hull $\XX (\vL)$
  has the desired properties. Indeed, the work of Robinson \cite{Rob}
  gives that, for any countable subgroup of $\RR^{d}$, one can find a
  cut and project scheme whose torus is just the dual of the
  subgroup. Then, any model set arising from a regular window from
  this cut and project scheme will be such a Meyer set
  \cite{Martin}. \exend
\end{remark}

It is possible to set up a diffraction theory for the elements of $\XX
(u)$ along the same lines as for $\XX (\vL)$. Indeed, if both $u$ and
$\vL$ are considered as translation bounded measures there is
virtually no difference in the framework and this is the point of view
proposed in \cite{BL-1}. As it is instructive, let us discuss the
diffraction theory of $u$. As before, we consider $u$ as a measure by
viewing it as a Radon--Nikodym density relative to Lebesgue measure
$\lambda$. Then, the measure $\widetilde{u \lambda}$ is given by
$\widetilde{u} \lambda$. Consequently, the autocorrelation of $u$ can
then simply be written as
\[
    \gamma^{}_{u} \, := \lim_{R\to \infty}
    \frac{u^{}_{R} \ast  \widetilde{u^{}_{R}}}
    {\vol (B^{}_{R} (0))} \ts ,
\]
where we use the shorthand $u^{}_{R} = u|^{}_{B_{R} (0)}$ for the
restriction of $u$ to the ball of radius $R$ around $0$. Of course,
the existence of the limit has still to be established.  Before we do
this, via an explicit calculation, let us pause for a very simple
special case.

\begin{example}\label{ex:one-diffraction}
  Consider $u\equiv 1$, hence Lebesgue measure itself. Then, a simple
  calculation with the volume-averaged convolution, compare
  \cite[Ex.~8.10]{TAO}, gives $\gamma^{}_{u} = \lambda$, and thus
  diffraction $\widehat{\gamma^{}_{u}} = \delta^{}_{0}$, which is a
  \emph{finite} pure point measure. Indeed, as we shall see later,
  this is an important distinction to the diffraction of a Delone set.
  \exend
\end{example}

To proceed with the general case, we will need two ingredients:
\begin{itemize}
\item One of the characteristic functions can be removed in the
  definition of $\gamma_u$. In particular, assuming existence of the
  limit, we have
\[
   \gamma^{}_{u}  \, := \lim_{R\to \infty}
   \frac{u^{}_{R} \ast \widetilde{u}}{\vol (B^{}_{R} (0))} \ts .
\]
(This is well-known and can be seen by a direct computation; compare
\cite{Martin,BL-1}).

\item For any $k\in\RR^{d}$, the limit
\[
   \lim_{R\to \infty} \frac{1}{\vol (B^{}_{R} (0))}
   \int_{B_{R} (0)}  \ee^{- 2 \pi \ii k x} u (x) \dd x
\]
exists. It is $a^{}_{k}$ if $k \in \cC$ and $0$ otherwise.  (This is
the formula for the Fourier--Bohr coefficient of $u$. It is easy to
see by direct computation and well-known in the theory of Bohr
almost periodic functions; see \cite{Cord,Katz} or
\cite[Thm.~8.2]{TAO}.)
\end{itemize}

Equipped with these two pieces of preparation, we are now going to
compute $\gamma^{}_{u}$. Let $g\in \cS$ be arbitrary. Using the first
ingredient, we find
\[
  \gamma^{}_{u} (g) \, =  \lim_{R\to \infty}
   \frac{\bigl( u^{}_{R} \nts \ast \widetilde{u} \ts \bigr) (g)}
   {\vol (B^{}_{R} (0))}  \ts .
\]
Direct computations then give
\begin{eqnarray*}
  \bigl( u^{}_{R} \nts \ast \widetilde{u}\ts  \bigr) (g)
  &=& \int_{\RR^{d}}  \int_{\RR^{d}}
  u^{}_{R} (y) \, \overline{u} (y-x) \, g (x) \dd y \dd x \\[2mm]
  &=& \int_{\RR^{d}}  \int_{\RR^{d}} u^{}_{R} (y)
   \sum_{k\in\cC} \overline{a^{}_{k}}\,
    \ee^{ - 2 \pi \ii  k (y-x)} g(x) \dd y \dd x  \\
   &=& \sum_{k\in\cC} \overline{a^{}_{k}}
    \int_{\RR^{d}} \int_{\RR^{d}} u^{}_{R} (y) \,
    \ee^{ - 2 \pi \ii k y}\ts \ee^{ 2 \pi \ii k x} g(x) \dd y \dd x \\
    &=&  \sum_{k\in\cC}
    \overline{a^{}_{k}} \int_{\RR^{d}}
    u^{}_{R} (y)\,  \ee^{ - 2 \pi \ii k y}
    \left( \int_{\RR^{d}} \ee^{ 2 \pi \ii k x} \ts g(x) \dd x \right) \dd y \\
   &=& \sum_{k\in\cC}  \overline{a^{}_{k}} \,
   F^{-1} (g)(k)  \int_{B_{R} (0)} u (y)
   \, \ee^{ -2 \pi \ii k y} \dd y \ts .
\end{eqnarray*}
Here, the second line follows from the definition of $u$, while
Fubini's theorem was employed in the penultimate step. Finally, last
step relies on the observation that the integral over $x$ just gives
the inverse Fourier transform $F^{-1} (g)$ of $g$.  Using the
preceding computation, the second ingredient and the summability of
the $(a^{}_{k})$, we then find
\[
   \gamma^{}_{u} (g) \, = \sum_{k\in \cC} |a^{}_{k} |^2
   \ts F^{-1} (g) (k) \ts .
\]
As this holds for all $g\in \cS$, we obtain
\[
   \gamma^{}_{u} \, = \, \biggl(\, \sum_{k\in \cC}
   |a^{}_{k} |^2 \,  \delta^{}_{k} \biggr) \circ {F^{-1} } \ts .
\]
Taking one more Fourier transform, and recalling $(\widehat{T}, g) =
(T, \widehat{g}\ts )$ for distributions $T$, we then find
\[
   \widehat{\gamma^{}_{u}} \, = \sum_{k\in \cC}
     |a^{}_{k} |^2 \, \delta^{}_{k} \ts .
\]
So, $\widehat{\gamma^{}_{u}}$ is a pure point measure with its set of
atoms being given by $\cC$.

\begin{remark}
  Due to the summability of the $(a^{}_{k})$, the $|a^{}_{k}|^2$ are
  also summable, and the pure point measure $\widehat{\gamma_u}$ is
  \emph{finite}, thus generalising the finding of
  Example~\ref{ex:one-diffraction}. In fact, one has the relation
\[
    \sum_{k \in \cC} \lvert a^{}_{k} \rvert^{2} \, =
    \lim_{R\to\infty} \frac{1}{\vol (B_{R} (0))}
    \int_{B_{R} (0)} \lvert u (x) \rvert^{2} \dd x \ts.
\]
This formula, which is not hard to derive from our above
considerations, is nothing but Parseval's identity for Bohr almost
periodic functions \cite[Thm.~I.1.18]{Cord}.  This way, one can see
immediately why the diffraction measure $\widehat{\gamma^{}_{u}}$ must
be a \emph{finite} measure.  This is an important structural
difference to the case of Delone sets. \exend
\end{remark}

Let us add the comment that this innocently looking observation, with
hindsight, sheds some light on the old dispute about the `right' model
for the description of quasicrystals between the quasiperiodic
function approach and the tiling or Delone set approach. While the
former leads to finite diffraction measures, the latter does not;
compare \cite[Rem.~9.11]{TAO} for a simple argument in the context of
cut and project sets, and \cite{Nicu+} as well as
\cite[Rem.~9.12]{TAO} for an argument in the more general situation of
Meyer sets.  Now, the experimental findings seem to indicate the
existence of series of Bragg peaks with growing $k$ and converging
intensity, which is not compatible with a finite diffraction measure
in the infinite volume limit.

\begin{remark}
  The Fourier--Bohr coefficients $a^{}_{k}$ as volume-averaged
  integrals can once again be interpreted as \emph{amplitudes} in our
  above sense, and it is then no surprise that the intensities of the
  Bragg peaks are once again given as the absolute squares of these
  amplitudes. This is another indication that there is more to be done
  in this direction.  \exend
\end{remark}

\section*{Acknowledgements}

The authors would like to thank the organizers of the \textit{3rd
  Bremen Winter School and Symposium:\ Diffusion on Fractals and
  Nonlinear Dynamics} (2015) for setting up a most stimulating event
inspiring in particular the material presented in
Section~\ref{Section:qpf}.  This work was supported by the German
Research Foundation (DFG), within the CRC 701.


\bigskip
\begin{thebibliography}{99}
\small

\bibitem{Aujogue}
J.-B.~Aujogue,
On embedding of repetitive Meyer multiple sets into
model multiple sets, to appear in:\
\textit{Ergodic Th.\ \& Dynam.\ Syst.};
\texttt{arXiv:1405.3013}.

\bibitem{ABKL}
J.-B.~Aujogue, M.~Barge, J.~Kellendonk and  D.~Lenz,
Equicontinuous factors, proximality and Ellis semigroup
for  Delone sets, in:\
J.~Kellendonk, D.~Lenz and J.~Savinien (eds): \textit{Mathematics of
Aperiodic Order}, Progress in Mathematics vol.\ 309, Birkh\"auser,
Basel (2015), pp.~137--194; \texttt{arXiv:1407.1787}.

\bibitem{Auslander88}
J.~Auslander,
\emph{Minimal Flows and their Extensions},
North Holland Mathematical Studies, vol.~153,
North Holland, Amsterdam (1988).

\bibitem{BG}
M.~Baake and F.~G\"{a}hler,
Pair correlations of aperiodic inflation rules via
renormalisation:\ Some interesting examples,
to appear in \textit{Top.\ Appl.};
\texttt{arXiv:1511.00885}.

\bibitem{BGG}
M.~Baake, F.~G\"{a}hler and U.~Grimm,
Examples of substitution systems and their factors
\textit{J.\ Int.\ Seq.} \textbf{16} (2013) art.\ 13.2.14 (18 pp);
\texttt{arXiv:1211.5466}.

\bibitem{TAO}
M.~Baake and U.~Grimm,
\textit{Aperiodic Order. Vol.\ $1$: A Mathematical Invitation},
Cambridge Univ.\ Press, Cambridge (2013).

\bibitem{squiral}
M.~Baake and U.~Grimm,
Squirals and beyond:\ Substitution tilings with singular
continuous spectrum,
\textit{Ergodic Th.\ \& Dynam.\ Syst.}
\textbf{34} (2014) 1077--1102;
\texttt{arXiv:1205.1384}.

\bibitem{BHS}
M.~Baake, C.~Huck and N.~Strungaru,
On weak model sets of extremal density,
\textit{preprint} \texttt{arXiv:1512.07129}.

\bibitem{BL-1}
M.~Baake and D.~Lenz,
Dynamical systems on translation bounded measures:\
Pure point dynamical and diffraction spectra,
\textit{Ergodic Th.\ \& Dynam.\ Syst.} \textbf{24} (2004)
1867--1893; \newline
\texttt{arXiv:math.DS/0302231}.

\bibitem{BL-2}
M.~Baake and D.~Lenz,
Deformation of Delone dynamical systems and
topological conjugacy,
\textit{J.\ Fourier Anal.\ Appl.} \textbf{11} (2005) 125--150;
\texttt{arXiv:math.DS/0404155}.

\bibitem{BLM}
M.~Baake, D.~Lenz and R.V.~Moody,
Characterization of model sets by dynamical systems,
\textit{Ergodic Th.\ \& Dynam.\ Syst.}
\textbf{27} (2007) 341--382;
\texttt{arXiv:math.DS/0511648}.

\bibitem{BLvE}
M.~Baake, D.~Lenz and A.C.D.~van Enter,
Dynamical versus diffraction spectrum for structures
with finite local complexity,
\textit{Ergodic Th.\ \& Dynam.\ Syst.}
\textbf{35} (2015) 2017--2043;
\texttt{arXiv:1307.7518}.

\bibitem{BM}
M.~Baake and R.V.~Moody,
Weighted Dirac combs with pure point diffraction,
\textit{J.\ reine angew.\ Math.\ (Crelle)}
\textbf{573} (2004) 61--94;
\texttt{arXiv:math.MG/0203030}.

\bibitem{BaakeMoody00}
M.~Baake and R.V.~Moody (eds),
\textit{Directions in Mathematical Quasicrystals},
CRM Monograph Series, vol.\ 13,
AMS, Providence, RI (2000).

\bibitem{BLPS}
S.~Beckus, D.~Lenz, F.~Pogorzelski and  M.~Schmidt,
Diffraction theory for processes of tempered distributions,
in preparation.

\bibitem{BF}
C.~Berg and G.~Forst,
\textit{Potential Theory on Locally Compact Abelian Groups},
Springer, Berlin (1975).

\bibitem{BT1}
E.~Bombieri and J.E.~Taylor,
Which distributions of matter diffract?
An initial investigation,
\textit{J.\ Phys.\ Colloques}
\textbf{47} (1986) C3-19--C3-28.

\bibitem{BT2}
E.~Bombieri and J.E.~Taylor,
Quasicrystals, tilings and algebraic numbers,
\textit{Contemp.\ Math.} \textbf{64} (1987) 241--264.

\bibitem{BDM}
X.~Bressaud, F.~Durand and A.~Maass,
Necessary and sufficient conditions to be an eigenvalue
for linearly recurrent dynamical Cantor systems,
\textit{J.\ London Math.\ Soc.} \textbf{72} (2005) 799--816;
\newline  \texttt{arXiv:0801.4619}.

\bibitem{Cord}
C.~Corduneanu,
\textit{Almost Periodic Functions}, 2nd English ed.\
Chelsea, New York (1989).

\bibitem{CFS}
I.P.~Cornfeld, S.V.~Fomin and Ya.G.~Sinai,
\textit{Ergodic Theory},
Springer, New York (1982).

\bibitem{Cow}
J.M.~Cowley,
\textit{Diffraction Physics},
3rd ed., North Holland, Amsterdam (1995).

\bibitem{Daley}
D.J.~Daley and D.~Vere-Jones,
\emph{An Introduction to the Theory of Point Processes},
Springer, New York (1988).

\bibitem{DM}
X.~Deng and R.V.~Moody,
Dworkin's argument revisited:\
point processes, dynamics, diffraction, and correlations,
\textit{J.\ Geom.\  Phys.} \textbf{58} (2008) 506--541;
\texttt{arXiv:0712.3287}.

\bibitem{Dekking}
F.M.~Dekking,
The spectrum of dynamical systems arising from substitutions
of constant length,
\textit{Z.\ Wahrscheinlichkeitsth.\ verw.\ Geb.}
\textbf{41} (1978) 221--239.

\bibitem{Dwo}
S.~Dworkin,
Spectral theory and $X\!$-ray diffraction,
\textit{J.\ Math.\ Phys.} \textbf{34} (1993) 2965--2967.

\bibitem{EW}
M.~Einsiedler and T.~Ward,
\textit{Ergodic Theory with a View towards Number Theory},
GTM 259, Springer, London (2011).

\bibitem{vEM}
A.C.D.~van Enter and J.~Mi\c{e}kisz,
How should one define a (weak) crystal?,
\textit{J.\ Stat.\ Phys.} \textbf{66} (1992) 1147--1153.

\bibitem{NF1}
N.~Frank,
Substitution sequences in $Z^d$ with a nonsimple Lebesgue component
in the spectrum,
\textit{Ergodic Th.\ \& Dynam.\ Syst.} \textbf{23} (2003) 519--532.

\bibitem{NF2}
N.P.~Frank,
Multidimensional constant-length substitution sequences,
\textit{Top.\ Appl.} \textbf{152} (2005) 44--69.

\bibitem{FR}
D.~Frettl\"{o}h and C.~Richard,
Dynamical properties of almost repetitive Delone sets,
\textit{Discr.\ Cont.\ Dynam.\ Syst.\ A} \textbf{34} (2014) 53--556;
\texttt{arXiv:1210.2955}.

\bibitem{GK}
F.~G\"ahler and R.~Klitzing,
The diffraction pattern of self-similar tilings,
in:\ R.V.~Moody (ed.),
\textit{The Mathematics of Long-Range Aperiodic Order},
NATO ASI Ser.\ C 489, Kluwer, Dordrecht (1997),
pp.~141--174.

\bibitem{GLA}
J.~Gil de Lamadrid and L.N.~Argabright,
Almost periodic measures,
\textit{Memoirs Amer.\ Math.\ Soc.} \textbf{85} (1990) no.~428
(AMS, Providence, RI).

\bibitem{Gou}
J.-B.~Gou\'{e}r\'{e},
Diffraction et mesure de Palm des processus ponctuels
(Diffraction and Palm measure of point prosesses)
\textit{C.\ R.\ Math.\ Acad.\ Sci.\ Paris, Ser.\ I}
\textbf{336} (2003) 57--62; \newline
\texttt{arXiv:math.PR/0208064}.

\bibitem{Gou-1}
J.-B.~Gou\'{e}r\'{e},
Quasicrystals and almost periodicity,
\textit{Commun.\ Math.\ Phys.} \textbf{255} (2005) 651--681;
\texttt{arXiv:math-ph/0212012}.

\bibitem{HvN}
P.R.~Halmos and J.~von Neumann,
Operator methods in classical mechanics. II.
\textit{Ann.\  Math.} \textbf{43} (1944) 332--350.

\bibitem{Hel}
H.~Helson,
Cocycles on the circle,
\textit{J.\ Operator Th.} \textbf{16} (1986) 189--199.

\bibitem{Hof}
A.~Hof,
On diffraction by aperiodic structures,
\textit{Commun.\ Math.\ Phys.} \textbf{169} (1995) 25--43.

\bibitem{Katz}
Y.~Katznelson,
\textit{An Introduction to Harmonic Analysis},
3rd ed., Cambridge University Press, Cambridge (2004).

\bibitem{Kel}
J.~Kellendonk,
Topological Bragg peaks and how they characterise point sets,
\textit{Acta Phys.\ Pol.\ A} \textbf{126} (2014) 497--500;
\texttt{arxiv:1309.7632}.

\bibitem{KL}
J.~Kellendonk and  D.~Lenz,
Equicontinuous Delone dynamical systems,
\textit{Can.\ J.\ Math.} \textbf{65} (2013) 149--170;
\texttt{arXiv:1105.3855}.

\bibitem{KS}
J.~Kellendonk and  L.~Sadun,
Meyer sets, topological eigenvalues, and Cantor fiber bundles,
\textit{J.\ London Math.\ Soc.}  \textbf{89} (2013) 114--130;
\texttt{arXiv:1211.2250}.

\bibitem{Koop}
B.O.~Koopman,
Hamiltonian systems and transformations in Hilbert space,
\textit{Proc.\ Nat.\ Acad.\ Sci.\ USA} \textbf{17} (1931) 315--318.

\bibitem{Lag-Delone}
J.C.~Lagarias,
Geometric models for quasicrystals I.~Delone sets of finite type,
\textit{Discr.\ Comput.\ Geom.} \textbf{21} (1999) 161--191.

\bibitem{Lag}
J.C.~Lagarias,
Mathematical quasicrystals and the problem of diffraction,
in:\ \cite{BaakeMoody00}, pp.~61--93.

\bibitem{LMS}
J.-Y.~Lee, R.V.~Moody and B.~Solomyak,
Pure point dynamical and diffraction spectra,
\textit{Ann.\ Henri Poincar\'e} \textbf{3} (2002) 1003--1018;
\texttt{arXiv:0910.4809}

\bibitem{LMS2}
J.-Y.~Lee, R.V.~Moody and B.~Solomyak,
Consequences of pure point diffraction spectra for
multiset substitution systems,
\textit{Discr.\ Comput.\ Geom.} \textbf{29} (2003) 525--560;
\texttt{arXiv:0910.4450}

\bibitem{Lenz}
D.~Lenz,
Continuity of eigenfunctions of uniquely ergodic
dynamical systems and intensity of Bragg peaks,
\textit{Commun.\ Math.\ Phys.} \textbf{287} (2009) 225--258;
\texttt{arXiv:math-ph/0608026}.

\bibitem{LM}
D.~Lenz and R.V.~Moody, Diffraction theory for stochastic processes,
\textit{preprint}, \texttt{arXiv:1111.3617 }


\bibitem{LR}
D.~Lenz and C.~Richard,
Pure point diffraction and cut and project schemes for
measures:\ The smooth case,
\textit{Math.\ Z.} \textbf{256} (2007) 347--378;
\texttt{math.DS/0603453}.

\bibitem{LSto}
D.~Lenz and  P.~Stollmann,
Delone dynamical systems and associated random operators,
in:\ J.-M.\ Combes, J.\ Cuntz, G.A.\ Elliott, G.\
Nenciu, H.\ Siedentop and S.\ Stratila (eds.),
\textit{Operator Algebras and Mathematical Physics},
Theta, Bucharest, (2003), pp.~267--285;
\texttt{arXiv:math-ph/0202042}.

\bibitem{LS}
D.~Lenz and N.~Strungaru,
Pure point spectrum for measure dynamical systems on locally
compact Abelian groups,
\textit{J.\ Math.\ Pures Appl.} \textbf{92} (2009) 323--341;
\texttt{arXiv:0704.2498}.

\bibitem{Loomis}
L.H.~Loomis,
\textit{An Introduction to Abstract Harmonic Analysis},
van Nostrand, Princeton (1953);
reprint, Dover, New York (2011).

\bibitem{M1}
Y.~Meyer,
\textit{Nombres de Pisot, Nombres de Salem et Analyse Harmonique},
LNM 117, Springer, Berlin (1970).

\bibitem{M2}
Y.~Meyer,
\textit{Algebraic Number Theory and Harmonic Analysis},
North Holland, Amsterdam (1972).

\bibitem{M3}
Y.~Meyer,
Quasicrystals, almost periodic patterns, mean-periodic functions
and irregular sampling,
\textit{African Diaspora J.\ Math.} \textbf{13} (2012) 1--45.

\bibitem{M-new}
Y.~Meyer,
Measures with locally finite support and spectrum,
\textit{preprint} (2016), to appear in \textit{PNAS}.

\bibitem{M-Nato}
R.V.~Moody,
Meyer sets and their duals,
in:\ R.V.~Moody (ed.),
 \textit{The Mathematics of Long-Range Aperiodic Order},
NATO ASI Series C 489, Kluwer, Dordrecht (1997), pp.~403--441.

\bibitem{M-beyond}
R.V.~Moody,
Model sets:\ A survey,
in:\ F.~Axel, F.~D\'{e}noyer and J.P.~Gazeau (eds.)
\textit{From Quasicrystals to More Complex Systems},
Springer, Berlin and EDP Sciences, Les Ulis (2000),
pp.~145--166; \texttt{arXiv:math.MG/0002020}.

\bibitem{MR}
P.~M\"{u}ller and C.~Richard,
Ergodic properties of randomly coloured point sets,
\textit{Can.\ J.\ Math.} \textbf{65} (2013) 349--402;
\texttt{arXiv:1005.4884}.

\bibitem{Phelps}
R.P.~Phelps,
\textit{Lectures on Choquet's Theorem},
2nd ed., LMN 1757, Springer, Berlin (2001).

\bibitem{Q}
M.~Queff\'{e}lec,
\textit{Substitution Dynamical Systems -- Spectral Analysis},
2nd ed., LNM 1294, Springer, Berlin (2010).

\bibitem{Rob}
E.A.~Robinson,
On uniform convergence in the Wiener--Wintner theorem,
\textit{J.\ London Math.\ Soc.} \textbf{49} (1994) 493--501.

\bibitem{Robbie}
E.A.~Robinson,
Symbolic dynamics and tilings of $\RR^{d}$,
\textit{Proc.\ Sympos.\ Appl.\ Math.} \textbf{60} (2004) 81--119.

\bibitem{Robbie2}
E.A.~Robinson,
A Halmos-von Neumann theorem for model set dynamical systems,
in:\ B.~Hasselblatt (ed.),
\textit{Dynamics, Ergodic Theory, and Geometry},
MSRI Publications, vl.~54, (2007), pp.~243--272.

\bibitem{Rudin}
W.~Rudin,
\textit{Fourier Analysis on Groups},
Wiley, New York (1962).

\bibitem{Danny}
D.~Shechtman, I.~Blech, D.~Gratias and J.W.~Cahn,
Metallic phase with long-range orientational order and
no translational symmetry,
\textit{Phys.\ Rev.\ Lett.} \textbf{53} (1984) 1951--1953.

\bibitem{Martin}
M.~Schlottmann,
Generalised model sets and dynamical systems,
in:\ \cite{BaakeMoody00}, pp.\ 143--159.

\bibitem{Schwartz}
L.~Schwartz,
\textit{Th\'{e}orie des Distributions},
reprint, Herman, Paris (1998).

\bibitem{Sol}
B.~Solomyak,
Dynamics of self-similar tilings,
\textit{Ergod.\ Th.\ \& Dynam.\ Syst.} \textbf{17} (1997), 695--738
and \textit{Ergod.\ Th.\ \& Dynam.\ Syst.} \textbf{19} (1999),
1685 (erratum).

\bibitem{Nicu}
N.~Strungaru, On weighted Dirac combs supported inside model sets,
\textit{J. Phys. A} \textbf{47}  (2014)  335202 (19 pp);
\texttt{arXiv:1309.7947}.

\bibitem{Nicu+}
N.~Strungaru,
On the Bragg diffraction spectra of a Meyer set,
\textit{Can.\ J.\ Math.} \textbf{65} (2013) 675--701;
\texttt{arXiv:1003.3019}.

\bibitem{ST}
N.~Strungaru and V.~Terauds,
Diffraction theory and almost periodic distributions,
\textit{preprint} (2016).

\bibitem{vN}
J.~von Neumann,
Zur Operatorenmethode in der klassischen Mechanik,
\textit{Ann.\ Math.} \textbf{33} (1933) 587--642.

\bibitem{Wal}
W.~Walter,
\textit{Einf\"{u}hrung in die Theorie der Distributionen},
3rd ed., BI-Wissenschaftsverlag, Mannheim (1994).

\end{thebibliography}
\end{document}